\documentclass[final]{siamltex}
\usepackage{xcolor}
\usepackage{amsmath}
\usepackage{amssymb}
\usepackage{mathrsfs}
\usepackage{graphicx}
\usepackage{bm}
\usepackage{ucs}
\usepackage[utf8x]{inputenc}
\usepackage{wrapfig}
\usepackage{color}
\usepackage{float}
\usepackage{showlabels}
\usepackage{epstopdf}

\usepackage[breaklinks,colorlinks=true,linkcolor=blue,citecolor=red,
  backref=page]{hyperref}

\newtheorem{algorithm}[theorem]{Algorithm}
\newtheorem{assumption}{Assumption}
\DeclareMathAlphabet{\itbf}{OML}{cmm}{b}{it}
\DeclareMathAlphabet\mathbfcal{OMS}{cmsy}{b}{n}

\def\RR{\mathbb{R}}
\def\CC{\mathbb{C}}
\def\PP{\mathbb{P}}
\def\bv{{{\itbf v}}}
\def\bx{{{\itbf x}}}
\def\vx{\vec{\bx}}
\def\bxp{\bx^\perp}
\def\by{{{\itbf y}}}
\def\vy{\vec{\by}}
\def\byp{\by^\perp}
\def\bz{{{\itbf z}}}
\def\vz{\vec{\bz}}
\def\bzp{\bz^\perp}

\newcommand{\vnu}{\vec{\boldsymbol{\nu}}_{_{\vx}}}

\newcommand{\cA}{{A}}
\newcommand{\cT}{\mathcal{T}}
\newcommand{\cM}{\mathcal{M}}
\newcommand{\cN}{\mathcal{N}}
\newcommand{\cS}{\mathcal{S}}
\newcommand{\cX}{\mathcal{X}}
\newcommand{\cW}{{W}}
\newcommand{\cF}{\mathcal{F}}
\newcommand{\csP}{{\mathscr P}}
\newcommand{\cK}{\mathcal{K}}
\newcommand{\us}{u^{\rm sc}}
\newcommand{\UsP}{U^{{\rm sc},\csP}}

\renewcommand{\theequation}{\thesection.\arabic{equation}}

\begin{document}


\title{Factorization method versus migration imaging in a
  waveguide} \author{Liliana Borcea and Shixu
  Meng\footnotemark[1]}

\maketitle

\renewcommand{\thefootnote}{\fnsymbol{footnote}}
\footnotetext[1]{Department of Mathematics, University of Michigan,
  Ann Arbor, MI 48109. {\tt borcea@umich.edu} and {\tt shixumen@umich.edu}}

\begin{abstract}
We present a comparative study of two qualitative imaging methods in
an acoustic waveguide with sound hard walls. The waveguide terminates
at one end and contains unknown obstacles of compact support, to be
determined from data gathered by an array of sensors that probe the
obstacles with waves and measure the scattered response. The first
imaging method, known as the factorization method, is based on the
factorization of the far field operator. It is designed to image at
single frequency and estimates the support of the obstacles by a Picard
range criterion.  The second imaging method, known as migration, works
either with one or multiple frequencies. It forms an
image by backpropagating the measured scattered wave to the search
points, using the Green's function in the empty waveguide.  
We study the connection between these methods  with analysis and numerical simulations.
\end{abstract}

\section{Introduction}
Qualitative approaches to inverse scattering problems have been the
focus of much activity in the mathematics community
\cite{cakoni2016qualitative,ammari2013mathematical}. Examples are the
linear sampling method
\cite{colton1996simple,arens2003linear,colton2003linear}, the
factorization method
\cite{kirsch1998characterization,kirsch2008factorization}, the
orthogonality sampling method
\cite{griesmaier2011multi,potthast2010study}, the range test method
\cite{potthast2003range}, and so on. Some of these methods are
connected to MUSIC (MUltiple-SIgnal-Classification)
\cite{cheney2001linear,kirsch2002music}, which is another qualitative
method that originates from signal processing
\cite{therrien1992discrete} and is used mostly for imaging point
scatterers
\cite{gruber2004time,borcea2016robust,moscoso2018robust,ammari2012direct}.

Reverse time migration methods and the closely related matched field
or matched filtering array data processing techniques are popular
in geophysics \cite{claerbout1985imaging,biondi}, ocean acoustics
\cite{bucker1976use,baggeroer1993overview}, radar imaging
\cite{Curlander,cheney2009fundamentals} and elsewhere. These methods
form an image by projecting data collected by a sensor array to the
replica wave field calculated for a point scatterer at the imaging
point. This projection is often called backpropagation. The high
frequency versions of these methods are based on the geometrical
optics approximation of the replica wave. They are known as Kirchhoff
migration \cite{biondi,bleistein2013mathematics} in broadband and
phase conjugation at a single frequency. 

Only some of the qualitative imaging methods, like orthogonality sampling 
\cite{griesmaier2011multi,potthast2010study}, are obviously related to
migration. The connection to the factorization method has been made
recently in \cite{liu2017novel}, for imaging in free space, using all
around measurements. Our goal in this paper is to extend these results
to imaging in a waveguide. 

Sensor array imaging in waveguides has  applications in
underwater acoustics \cite{baggeroer1993overview}, imaging of and in
tunnels \cite{schultz2007remote,haack1995state,bedford2014modeling},
nondestructive evaluation of slender structures
\cite{rizzo2010ultrasonic}, and so on. Migration type imaging methods
in waveguides with perfectly known geometry have been developed and
analyzed in
\cite{dediu2006recovering,bourgeois2008linear,monk2012sampling,
  monk2016inverse,borcea2015imaging,tsogka2013selective,tsogka2017imaging}
and examples of imaging with experimental validation are in
\cite{mordant1999highly,philippe2008characterization}. The case of
unknown waveguide geometry is more difficult and is addressed in
\cite{borcea2013quantitative,borcea2018ghost} for randomly perturbed
waveguide boundary. We also refer to \cite{borcea2018direct} for a
linear sampling approach to imaging in a waveguide with unknown,
compactly supported wall deformations. Linear sampling imaging in
waveguides with known geometry is studied in
\cite{xu2000generalized,bourgeois2008linear,bourgeois2011use,
bourgeois2012use,monk2012sampling}.

We are interested in the factorization method and its connection to
migration, for imaging obstacles in a waveguide with known geometry,
that terminates at one end. The termination is motivated by the
application of imaging in tunnels and is beneficial because the
reflection at the end wall allows a back view of the obstacles. The main
difference between the factorization method in a waveguide and in free
space is due to the fact that in the waveguide the wave field is a
superposition of finitely many propagating modes and infinitely many
evanescent modes which cannot be measured in the far field. Thus,
imaging must be done only with the propagating modes.  

So far, the factorization method in waveguides and cavities has been restricted to
using unphysical incident waves as explained in \cite[Section
  1.7]{kirsch2008factorization} and
\cite{arens2011direct,bourgeois2014identification,meng2014factorization}.
This issue is addressed in \cite{bourgeois2012use}, by considering
incident fields that are pure guided modes and measuring the reflected
and transmitted modes before and after the obstacle. Such incident
fields could be obtained with a full aperture array of sources, but
the measurement of the reflected and transmitted modes may be
difficult to realize in some applications.

In this paper we show that the factorization method can be used in a
terminated waveguide, for physical incident waves generated by sensors
in an array that lies far from the obstacle, on the opposite side of
the end wall. We establish a connection between the factorization method 
and migration imaging and show that obstacles can be localized using 
only the propagating part of the wave field.

The paper is organized as follows: We begin in section
\ref{sect:formulation} with the formulation of the inverse scattering
problem. Then, we discuss in section \ref{sect:factorization} the
factorization method. The connection to migration imaging is in
section \ref{sect:migration}. We assess the results with numerical
simulations in section \ref{sect:numerics} and end with a summary
in section \ref{sect:summary}.

\section{The inverse problem}
\label{sect:formulation}

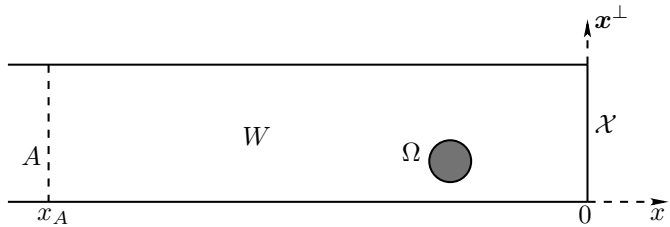
\begin{figure}[t]
\begin{center}
\includegraphics[width=0.68\textwidth]{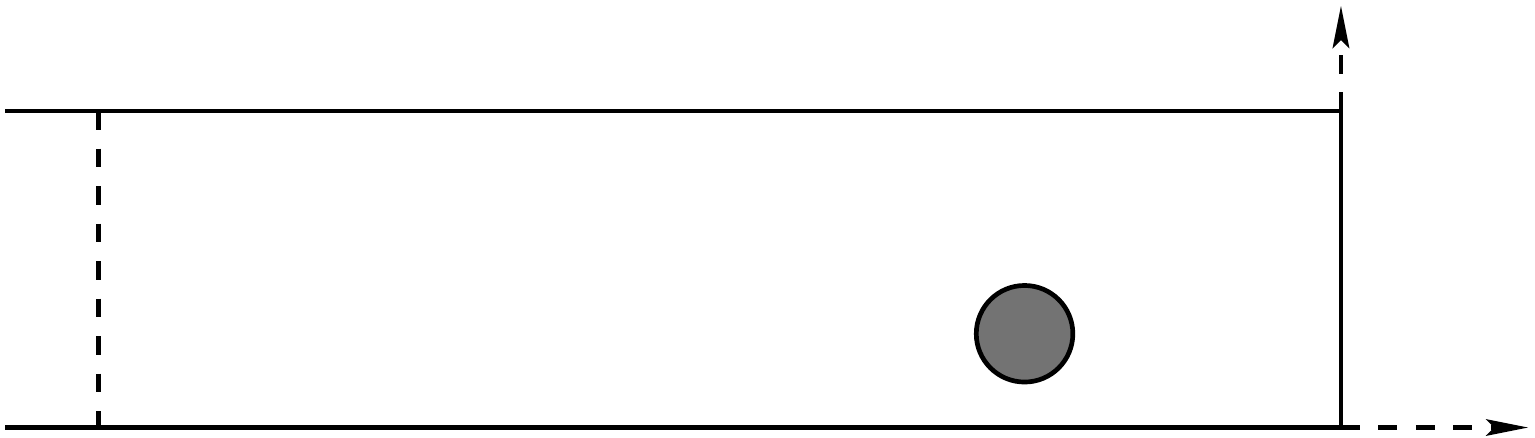}%
\end{center}
\setlength{\unitlength}{3947sp}%
\begingroup\makeatletter\ifx\SetFigFont\undefined%
\gdef\SetFigFont#1#2#3#4#5{%
  \reset@font\fontsize{#1}{#2pt}%
  \fontfamily{#3}\fontseries{#4}\fontshape{#5}%
  \selectfont}%
\fi\endgroup%
\begin{picture}(9591,2871)(1918,-3049)
\put(3000,100){\makebox(0,0)[lb]{\smash{{\SetFigFont{7}{8.4}{\familydefault}{\mddefault}{\updefault}{\color[rgb]{0,0,0}{\normalsize $\cA$ }}%
}}}}
\put(3100,-250){\makebox(0,0)[lb]{\smash{{\SetFigFont{7}{8.4}{\familydefault}{\mddefault}{\updefault}{\color[rgb]{0,0,0}{\normalsize $x_{\cA}$ }}%
}}}}
\put(6950,-250){\makebox(0,0)[lb]{\smash{{\SetFigFont{7}{8.4}{\familydefault}{\mddefault}{\updefault}{\color[rgb]{0,0,0}{\normalsize $x$ }}%
}}}}
\put(6600,960){\makebox(0,0)[lb]{\smash{{\SetFigFont{7}{8.4}{\familydefault}{\mddefault}{\updefault}{\color[rgb]{0,0,0}{\normalsize $\bxp$ }}%
}}}}
\put(6500,-270){\makebox(0,0)[lb]{\smash{{\SetFigFont{7}{8.4}{\familydefault}{\mddefault}{\updefault}{\color[rgb]{0,0,0}{\normalsize $0$ }}%
}}}}
\put(5390,120){\makebox(0,0)[lb]{\smash{{\SetFigFont{7}{8.4}{\familydefault}{\mddefault}{\updefault}{\color[rgb]{0,0,0}{\normalsize $\Omega$ }}%
}}}}
\put(4390,220){\makebox(0,0)[lb]{\smash{{\SetFigFont{7}{8.4}{\familydefault}{\mddefault}{\updefault}{\color[rgb]{0,0,0}{\normalsize $\cW$ }}%
}}}}
\put(6600,300){\makebox(0,0)[lb]{\smash{{\SetFigFont{7}{8.4}{\familydefault}{\mddefault}{\updefault}{\color[rgb]{0,0,0}{\normalsize $\cX$ }}%
}}}}
\end{picture}%
\vspace{-2.4in}
\caption{Imaging setup: An obstacle supported in $\Omega$ in the
  waveguide $\cW = (-\infty, 0) \times \cX$ is imaged using
  measurements collected by an array of sensors lying in the set
  $\cA$, at range offset $|x_{\cA}|$ from the end wall. The system of
  coordinates $\vx = (x,\bxp)$ is centered at the end wall, with range
  coordinate $x<0$ in the waveguide $\cW$ and cross-range coordinate
  $\bxp$ in the cross-section $\cX$.}
\label{fig:setup}
\end{figure}

Consider a waveguide that terminates at one end
\begin{equation}
\label{eq:WG}
\cW = (-\infty,0) \times \cX \subset \RR^{d}, \qquad 2 \le d \le 3,
\end{equation}
with cross-section $\cX \subset \RR^{d-1}$. In two dimensions $(d=2)$
$\cX$ is the interval $(0,|\cX|)$ of length $|\cX|$, whereas in three
dimensions $\cX$ is a convex and bounded domain with piecewise smooth boundary
$\partial \cX$. We use  the system of coordinates $\vx =
(x,\bxp)$ with range $x$ along the axis of the waveguide,
starting from the end wall, and with cross-range $\bxp \in \cX$.  To
fix ideas, we assume that the waveguide has sound hard walls 
\begin{equation}
\label{eq:WGwall}
\partial \cW = \{0\} \times \cX \cup (-\infty,0) \times \partial \cX,
\end{equation}
and contains sound soft obstacles supported in the compact set
$\Omega \subset \cW$, with piecewise smooth boundary $\partial
\Omega$. The results are expected to extend to other boundary conditions at $\partial
\cW$ and $\partial \Omega$, and also to penetrable scatterers.

The inverse scattering problem is to determine the obstacles from
measurements gathered by an array of $n_{\cA}$ sensors located in the set
\begin{equation}
\label{eq:Array}
\cA = \{x_{\cA} \} \times \cX, \qquad x_{\cA} < 0,
\end{equation}
that lies on the left side of the obstacles, as illustrated in Figure
\ref{fig:setup}. For simplicity of the presentation we carry out the
analysis in the full aperture case\footnote{The factorization method with a partial aperture array requires additional data processing,
as explained in section \ref{sect:numerics} and in \cite[Section 2.4]{borcea2018direct}, whereas the  implementation of the migration method is independent of the aperture. }, where the array spans
the entire set $\cA$.   

The array probes the waveguide with a time harmonic wave emitted from
one of the sensors, at location $\vx_s \in \cA$, and measures the
echoes $\us(\vx_r,\vx_s)$ at all the sensor locations $\vx_r \in \cA$.
These echoes are defined in section \ref{sect:form.2}. The array data is the 
 response matrix
\begin{equation}
\label{eq:RESPM}
\itbf{U}^{\rm sc} = \big(\us(\vx_r,\vx_s)\big)_{1 \le r,s \le n_{\cA}},
\end{equation}
gathered by successive illuminations, with one source at a time.
We assume in the analysis that the sensor spacing is sufficiently
small, so we can make the continuum aperture approximation. This means
that  we replace sums over the source and receiver
indexes $s,r = 1, \ldots, n_{\cA}$ by integrals over the aperture $\cA$.

\subsection{The incident wave}
\label{sect:form.1}
The probing (incident) wave emitted by the source at $\vx_s \in \cA$
is defined by the solution of the Helmholtz equation in the empty
waveguide. It is the Green's function $G(\vx,\vx_s)$ satisfying
\begin{align}
\big(\Delta_{\vx} + k^2\big) G(\vx,\vx_s) &= - \delta(\vx-\vx_s),  \qquad 
\vx \in \cW, \nonumber \\
\partial_{\vnu}G(\vx,\vx_s) &= 0, \hspace{0.9in}\vx \in \partial \cW, 
\label{eq:GrF}
\end{align}
and the outgoing radiation condition at range $x < x_{\cA}$, stated in
Definition \ref{def.1}. Here $\Delta_{\vx}$ is the Laplace operator,
$k$ is the wavenumber and $\vnu$ denotes the normal at $\partial \cW$
at point $\vx \in \partial \cW$.

\vspace{0.05in}
\begin{definition}
\label{def.1}
We say that a time harmonic wave field $v(\vx)\exp(-i \omega t)$,
where $\omega$ is the frequency and $t$ is time, satisfies the
``outgoing radiation condition'' at range $x$ if it consists of
backward (left) going modes and decaying evanescent modes. The wave
satisfies the ``incoming radiation condition'' at range $x$ if it
consists of forward (right) going modes and decaying evanescent modes.
\end{definition}

\vspace{0.05in}
The mode decomposition of the Green's function is obtained via
separation of variables i.e., by expansion in the $L^2(\cX)$
basis $\{\psi_j(\bxp)\}_{j \ge 0}$ of eigenfunctions of the Laplace operator
$\Delta_{\bxp}$ in the cross-range $\bxp$, with Neumann boundary
conditions at $\partial \cX$. 
These eigenfunctions can be chosen to be real-valued. They 
satisfy
\begin{align}
\Delta_{\bxp} \psi_j(\bxp) &= - \lambda_j \psi_j(\bxp), \qquad \bxp
\in \cX, \nonumber \\ \partial_{\boldsymbol{\nu}_{\bxp}} \psi_j(\bxp)
&= 0,  \hspace{0.9in}\bxp \in \partial \cX, 
\label{eq:EigF}
\end{align}
and  are orthonormal 
\begin{equation}
\int_{\cX} d \bxp \, \psi_j(\bxp) \psi_{j'}(\bxp) = \delta_{j,j'}.
\label{eq:EigFOrt}
\end{equation}
The eigenvalues $-\lambda_j$ are real and are ordered as $ 0 = \lambda_o <
\lambda_1 \le \lambda_2 \le \ldots. $ They determine the number $J+1$ of 
propagating modes, where 
\begin{equation}
J = \max\{j \in \mathbb{N}: ~ \lambda_j \le k^2\}.
\label{eq:J}
\end{equation}
The modes indexed by $j = 0, \ldots, J$ are one dimensional time
harmonic waves of the form $\exp[ i (\pm  \beta_j x - \omega t)]$
propagating forward (to the right) and backward (to the left) along
the range direction $x$, with wavenumber
\begin{equation}
\beta_j = \sqrt{k^2 - \lambda_j}, \qquad j = 0, \ldots, J.
\label{eq:betaP}
\end{equation}
The infinitely many modes indexed by $j > J$ are evanescent waves that
decay exponentially away from the source, on the range scale
$1/|\beta_j|$, where
\begin{equation}
\beta_j = 
i \sqrt{\lambda_j-k^2}, \qquad  j > J.
\label{eq:betaE}
\end{equation}

We assume throughout that the probing frequency is such that $\beta_j
\ne 0$ for all $j \ge 0$. Then, the  incident field due to the source at $\vx_s = (x_{\cA},\bxp_s) \in \cA$ is given by 
\begin{align}
u^{\rm inc}(\vx,\vx_s) = G(\vx,\vx_s) &= \sum_{j=0}^J \frac{i}{
  \beta_j} \psi_j(\bxp) \psi_j(\bxp_s) e^{- i \beta_j x_{\cA}}
\cos(\beta_j x)\nonumber \\ & + {\sum_{j>J}} \frac{1}{|\beta_j|}
\psi_j(\bxp) \psi_j(\bxp_s) e^{|\beta_j| x_{\cA}} \cosh(|\beta_j| x),
\label{eq:uinc}
\end{align}
at $\vx = (x,\bxp)$ on the right of the array, with range $x \in (x_{\cA},0)$.
At points on the left of the array, with  range $x < x_{\cA}$, the expression of
$G(\vx,\vx_s)$ is obtained by interchanging $x$ with $x_{\cA}$ in the
right hand side of \eqref{eq:uinc}.

Note that $u^{\rm inc}(\vx,\vx_s) \exp(-i \omega t)$ satisfies the
outgoing radiation condition at range $x < x_{\cA}$, whereas between
the array and the end wall there are both
forward and backward propagating modes. Because we assume a fixed
frequency $\omega$ in the analysis, we drop henceforth the factor
$\exp(-i \omega t)$.
\subsection{The scattered wave}
\label{sect:form.2}
To define the scattered wave, we make the following standard
assumption:

\vspace{0.05in}
\begin{assumption}
\label{as.1}
The wavenumber $k$ is such that the problem 
\begin{align*}
(\Delta_{\vx} + k^2) w(\vx) &=0, \qquad \vx \in \cW \setminus
  \overline{\Omega}, \\ \partial_{\vnu} w(\vx) &= 0, \hspace{0.29in} \vx
  \in \partial \cW, \\ w(\vx)&=0, \hspace{0.29in} \vx \in \partial
  \Omega,
\end{align*}
has only the trivial solution $w(\vx) \equiv 0$ that satisfies either
the outgoing or the incoming radiation condition on the left side of
$\Omega$, at range 
\[
x < x_{\Omega} = \inf\{x: \vx = (x,\bxp) \in \Omega\}.
\]
Here $\overline{\Omega}$ denotes the closure of $\Omega$.
\end{assumption}

\vspace{0.05in}
With this assumption, it is known (see for example \cite[Theorem
  A.4]{borcea2018direct}) that the scattered wave field
$\us(\vx,\vx_s)$, satisfying
\begin{align}
(\Delta_{\vx} + k^2) \us(\vx,\vx_s) &=0, \qquad \qquad ~~~\vx \in \cW
  \setminus \overline{\Omega}, \label{eq:usc1}\\ \partial_{\vnu}
  \us(\vx,\vx_s) &= 0, \hspace{0.69in} \vx \in \partial
  \cW,\label{eq:usc2} \\ \us(\vx,\vx_s)&=-G(\vx,\vx_s), ~~~ \vx \in
  \partial \Omega, \label{eq:usc3}
\end{align}
and the outgoing radiation condition at range $x < x_{\Omega}$, is
well defined. Moreover, $\us(\cdot, \vx_s) \in H_{\rm
  loc}^1(\cW\setminus \overline{\Omega})$.

We will need a second assumption, which holds for
all positive $k$ with the exception of a countable set:

\vspace{0.05in}
\begin{assumption}
\label{as.2}
The wavenumber $k$ is such $k^2$ is not an eigenvalue of the
negative Laplacian in $\Omega$ with Dirichlet boundary conditions at $\partial \Omega$. That is to say, the problem
\begin{align*}
(\Delta_{\vx} + k^2) w(\vx) &=0, \qquad \vx \in \Omega, \\
w(\vx)&=0, \hspace{0.29in} \vx \in \partial \Omega,
\end{align*}
has only the trivial solution $w(\vx) \equiv 0$ in $\Omega$.
\end{assumption}

\section{Imaging with the factorization method}
\label{sect:factorization}
We now describe  the factorization method for solving the inverse
scattering problem.  We begin
in section \ref{sect:fact.1} with the definition of the relevant
operators and then describe the method in section \ref{sect:fact.2}.

\subsection{The operators}
\label{sect:fact.1}
Consider the linear integral operator $\cN:L^2(\cA) \to L^2(\cA)$,
\begin{equation}
\cN g(\vx) = \int_{\cA} d S_{\vy} \, \us(\vx,\vy) g(\vy), \qquad
\vx \in \cA, ~~\forall \, g \in L^2(\cA),
\label{eq:F1}
\end{equation}
with kernel given by the measured scattered field $\us$ at the array. This is called in
the literature, depending on the authors, either the far field or the
near field operator. It defines the scattered wave $\cN g$
received at the array, due to an illumination $g$ from all the sources
in $\cA$.  Because $\us(\cdot, \vx_s) \in H_{\rm loc}^1(\cW\setminus
\overline{\Omega})$, the range of $\cN$ lies in
$H^{\frac{1}{2}}(\cA)$, but we view $\cN$ as an operator from
$L^2(\cA)$ to $L^2(\cA)$.  As shown in the next section, $\cN$ can be
factorized in terms of three linear operators $\cT$, $\Lambda$ and $\cS$ that we now define:

The operator $\cT:L^2(\cA) \to H^{\frac{1}{2}}(\partial \Omega)$
maps functions defined at the array  to functions
defined at the boundary $\partial \Omega$ of the obstacles,
\begin{equation}
\cT g(\vz) = \int_{\cA} d S_{\vy} \, G(\vz,\vy) g(\vy), \qquad
 \vz \in \partial \Omega, ~~ \forall \, g \in L^2(\cA).
\label{eq:F2}
\end{equation}
Its adjoint $\cT^\star:H^{-\frac{1}{2}}(\partial \Omega)\to L^2(\cA)$ is
given by
\begin{equation}
\cT^\star h(\vx) = \int_{\partial \Omega} d S_{\vz} \,
\overline{G(\vz,\vx)} h(\vz), \qquad \vx \in \cA, ~~
\forall \, h \in H^{-\frac{1}{2}}(\partial \Omega),
\label{eq:F3}
\end{equation}
where the bar denotes throughout the complex conjugate. This adjoint is defined using
the inner product
\begin{equation}
\label{eq:innprodA}
\big( f,g \big)_{\cA} = \int_{\cA} dS_{\vx} \, \overline{f(\vx)}
g(\vx), \qquad \forall \, f,g \in L^2(\cA),
\end{equation}
and the duality pairing 
\begin{equation}
\label{eq:dualityOm}
\left< f,g\right>_{\partial \Omega} = \int_{\partial \Omega} dS_{\vx}
\, \overline{f(\vx)} g(\vx), \qquad \forall \, f \in
H^{\frac{1}{2}}(\partial \Omega), ~~ \forall \, g \in H^{-\frac{1}{2}}(\partial
\Omega),
\end{equation}
meaning that
\begin{equation}
\label{eq:dualT}
\big( f, \cT^\star g\big)_{\cA} = \left< \cT f, g \right>_{\partial
  \Omega}, \qquad \forall \, f \in L^2(\cA), ~~\forall \,  g \in
H^{-\frac{1}{2}}(\partial \Omega).
\end{equation}

The operator  $\Lambda:
H^{\frac{1}{2}}(\partial \Omega) \to  H^{-\frac{1}{2}}(\partial
\Omega)$ is the Dirichlet to Neumann map 
\begin{equation}
\Lambda f(\vx) = h(\vx), \qquad  \vx \in \partial \Omega, ~ ~
\forall \, f \in H^{\frac{1}{2}}(\partial \Omega),
\label{eq:F4}
\end{equation}
where $h\in H^{-\frac{1}{2}}(\partial \Omega)$ is the solution of 
\begin{equation}
\int_{\partial \Omega} dS_{\vz} \, \overline{G(\vx,\vz)} h(\vz) = -
f(\vx), \qquad  \vx \in \partial \Omega.
\label{eq:F5}
\end{equation}
The solvability of \eqref{eq:F5} is established in \cite[Section
  4.2]{bourgeois2008linear}, under the Assumption \ref{as.2}, 
and \cite[Proposition 1]{bourgeois2008linear} gives that 
$\Lambda$ is an isomorphism.

The  scattering operator $\cS$ maps incoming to outgoing waves at $\cA$. To define it, we
introduce the function spaces
\begin{align}
\mathscr{H}(W\setminus \overline{\Omega}) &= \Big\{ w \in H_{\rm
  loc}^1(W\setminus \overline{\Omega}): ~ ~ (\Delta_{\vx} + k^2)
w(\vx) = 0 ~{\rm in } ~ W\setminus \overline{\Omega}, \nonumber
\\&\hspace{1.9in}  \partial_{\vnu} w(\vx) = 0 ~{\rm on } ~
\partial W \Big\},
\label{eq:F6} \\
\mathscr{H}^{\rm out}(A) &= \Big\{ w|_{\cA}: ~ w \in
\mathscr{H}(W\setminus \overline{\Omega}), ~~ w \mbox{
  satisfies the outgoing } \nonumber \\ &\hspace{1.4in}\mbox{radiation
  condition at } x < x_\Omega \Big\},
\label{eq:F7}
\\ \mathscr{H}^{\rm in}(A) &= \Big\{ w|_{\cA}: ~ w \in
\mathscr{H}(W\setminus \overline{\Omega}), ~~ w \mbox{
  satisfies the incoming} \nonumber \\ &\hspace{1.4in}\mbox{radiation
  condition at } x < x_\Omega \Big\},
\label{eq:F8}
\end{align}
where $w|_{\cA}$ denotes the trace of $w$ on $\cA$. The operator  $\cS :
\mathscr{H}^{\rm in}(A) \to \mathscr{H}^{\rm out}(A)$ is defined by
\begin{align}
\cS v(\vx) = w(\vx), \qquad  \vx \in \cA, ~ ~ \forall \, v
\in \mathscr{H}^{\rm in}(A),
\label{eq:F9}
\end{align}
where $w(\vx) \in \mathscr{H}(W\setminus \overline{\Omega})$ satisfies the 
boundary condition 
\begin{equation}
w(\vx) = v(\vx), \qquad \vx \in \partial \Omega,
\label{eq:F10}
\end{equation}
and the outgoing radiation condition at range $x <
x_{\Omega}$. Moreover,  $\cS$ is invertible\footnote{This
  follows by the unique solvability of the Helmholtz equation in $\cW
  \setminus \overline{\Omega}$ with homogeneous Neumann conditions at
  $\partial \cW$ and outgoing or incoming radiation condition, using
  that $v|_{\partial \Omega} = w_{\partial \Omega}$.}, with inverse
$\cS^{-1} : \mathscr{H}^{\rm out}(A) \to \mathscr{H}^{\rm in}(A)$
defined by
\begin{align}
\cS^{-1} w(\vx) = v(\vx), \qquad  \vx \in \cA, ~ ~ \forall
\, w \in \mathscr{H}^{\rm out}(A),
\label{eq:F11}
\end{align}
where $v(\vx) \in \mathscr{H}(W\setminus \overline{\Omega})$ satisfies the 
boundary condition 
\begin{equation}
{v(\vx) = w(\vx)}, \qquad \vx \in \partial \Omega,
\label{eq:F12}
\end{equation}
and the incoming radiation condition at range $x <
x_{\Omega}$.
\subsection{The factorization method}
\label{sect:fact.2}
The imaging  is based on the operator
\begin{equation}
\cF:L^2(\cA) \to L^2(\cA), ~ ~~ \cF = \cS^{-1} \cN,
\label{eq:F14}
\end{equation}
which is defined in terms of the array measurements, as stated in the
following lemma:
\begin{lemma}
\label{lem.1}
Any $\phi \in L^2(\cA)$ can be written as  
\begin{equation}
\phi(\vx) = \phi^{(1)}(\vx) + i \phi^{(2)}(\vx), 
\label{eq:F15}
\end{equation}
for $\vx = (x_{\cA}, \bxp) \in \cA$, with
\begin{equation}
\phi^{(l)}(\vx) = \sum_{j=0}^J \frac{\alpha_j^{(l)}}{i} \psi_j(\bxp)
e^{i \beta_j x_{\cA}} + \sum_{j > J} \alpha_j^{(l)} \psi_j(\bxp),
\label{eq:F16}
\end{equation}
and $\alpha_j^{(l)} \in \RR$, for all $j \ge 0$ and $l = 1,2$.  Furthermore,
\begin{equation}
\cF \phi(\vx) = \cF \phi^{(1)}(\vx) + i \cF \phi^{(2)}(\vx),
\label{eq:F17}
\end{equation}
where 
\begin{equation}
\cF \phi^{(l)}(\vx) = \int_{\cA} dS_{\vy} \, \overline{\us(\vx,\vy)
  \phi^{(l)}(\vy)}, \qquad \vx \in \cA, ~ ~ l = 1,2.
\label{eq:F18}
\end{equation}
\end{lemma}

The proof of this lemma is in Appendix \ref{ap:A} and the
decomposition \eqref{eq:F15} is obtained from the expansion of $\phi$
in the $L^2(\cA)$ eigenbasis $\{\psi_j(\bxp)\}_{j \ge 0},$
\begin{equation}
\phi(\vx) = \sum_{j=0}^\infty \gamma_j \psi_j(\bxp), \qquad 
\vx = (x_{\cA}, \bxp) \in \cA,
\label{eq:F19}
\end{equation}
with coefficients $\gamma_j \in \CC$. The real valued $\alpha_j^{(1)}$
and $\alpha_j^{(2)}$ in \eqref{eq:F16} are defined in terms of these
coefficients by
\begin{align}
\alpha_j^{(1)} + i \alpha_j^{(2)} &= \left\{ \begin{array}{ll} i 
  \gamma_j e^{-i \beta_j x_\cA}, \qquad &\mbox{if} ~ ~ j = 0, \ldots,
  J,\\ 
  \gamma_j, &  \mbox{if} ~~ j > J.
\end{array} \right. \label{eq:F21}
\end{align}

\vspace{0.05in}
\begin{theorem}
\label{thm.1}
The operator $\cF$ has the factorization
\begin{equation}
\cF =  \cT^\star \Lambda \cT,
\label{eq:F13}
\end{equation}
and the operators defined in (\ref{eq:F2}) and \eqref{eq:F4}) satisfy the
following properties: \\ (i) The operator $\cT$ is compact and 
injective.  \\(ii) Let $\Lambda^\star:H^{\frac{1}{2}}(\partial \Omega) \to
H^{-\frac{1}{2}}(\partial \Omega)$ be the adjoint of $\Lambda$, defined by 
\begin{equation}
\left< \Lambda f, g \right>_{\partial \Omega} = \left< f,
\Lambda^\star g \right>_{\partial \Omega}, \qquad \forall \, f,g \in
H^{\frac{1}{2}}(\partial \Omega),
\end{equation}
using the duality pairing \eqref{eq:dualityOm}.
Define the self-adjoint operators $ \Im(\Lambda) = \big(\Lambda-\Lambda^\star\big)/(2i)$ and 
$\Re(\Lambda) = \big(\Lambda+\Lambda^\star\big)/2$. Then, $-\Im(\Lambda)$ is positive semi-definite,
\begin{align}
{-\left< \Im(\Lambda) f, f \right>_{\partial \Omega}} \ge 0,
\qquad \forall \, f \in H^{\frac{1}{2}}(\partial \Omega), 
\label{eq:F22}
\end{align}
and  ${-\Re(\Lambda) }$
is the sum of a positive definite, self-adjoint operator and a compact
operator.
\end{theorem}

\vspace{0.05in}
This result, proved in Appendix \ref{ap:B}, and the next lemma, proved
in Appendix \ref{ap:C}, are the theoretical foundation of the
factorization method.

\vspace{0.05in}
\begin{lemma}
\label{lem.2}
Let $\vz \in (x_{\cA},0) \times \cX$ be a search point. Then, $\vz \in
\Omega$ if and only if $\overline{G(\cdot,\vz)}|_{\cA} \in
\mbox{range}(\cT^\star).$
\end{lemma}
\vspace{0.05in}

The range test in Lemma \ref{lem.2} cannot be used directly to determine the
support $\Omega$ of the obstacles, because $\cT^\star$ is not
known. However, \cite[Theorem 2.1]{lechleiter2009factorization} shows
that $\Omega$ can be determined using a new operator 
\begin{equation}
\cF_{\#} =  \big|\Re(\cF)\big| -\Im(\cF) 
\label{eq:F23} 
\end{equation}
where 
\begin{equation}
\Re(\cF) = \frac{\big(\cF + \cF^\star\big)}{2}, \qquad 
\Im(\cF) = \frac{\big(\cF - \cF^\star\big)}{2i},
\label{eq:F13p}
\end{equation}
and $\big|\Re(\cF)\big|$ is defined in the standard way,
using the spectral representation of $\Re(\cF)$. Similarly, we define 
\begin{align} 
\Re(\Lambda) = \frac{(\Lambda + \Lambda^\star)}{2}, \qquad \Im(\Lambda) = \frac{ (\Lambda - \Lambda^\star)}{2i}, \qquad 
\Lambda_\# = \big| \Re(\Lambda) \big| - \Im(\Lambda),
\label{eq:F13pp}
\end{align}
and conclude from the proof of \cite[Theorem 2.1]{lechleiter2009factorization} that 
\begin{equation}
\cF_{\#} =   \cT^* \Lambda_\# \cT.
\label{eq:F23-1} 
\end{equation}
We deduce from Theorem \ref{thm.1} and \eqref{eq:F23} that $\cF_{\#}$ is positive definite, so we can
take its square root $\cF_{\#}^{\frac{1}{2}}$.  The following result follows from  Theorem \ref{thm.1}, Lemma \ref{lem.2} and
\cite[Theorem 2.1]{lechleiter2009factorization}.

\vspace{0.05in}
\begin{theorem}
\label{thm.Fsharp}
Let $\vz \in (x_{\cA},0) \times \cX$ be a search point in the waveguide, between the array and the end wall. Then, $\vz \in \Omega$ if and 
only if 
\begin{equation}
\mbox{inf}\Big\{ \big( {  \cF^\#} \varphi,\varphi \big)_{\cA}:
~ \varphi \in L^2(\cA), ~ ~ \big(\overline{G(\cdot,
  \vz)},\varphi \big)_{\cA} = 1 \Big\} > 0,
\label{eq:FR6}
\end{equation}
or, equivalently, if and only if 
\begin{equation}
\overline{G(\cdot,\vz)}|_{\cA} \in
\mbox{range}\big(\cF_{\#}^{\frac{1}{2}}\big).
\label{eq:F24}
\end{equation}
\end{theorem}

The factorization method uses  the condition \eqref{eq:F24} and a Picard range criterion
to define the sampling function
\begin{equation}
g_{\#}(\vz) = \sum_{j=1}^\infty \frac{\big|
\big( \overline{G(\cdot,\vz)}, \varphi_j\big)_{\cA}\big|^2}{\mu_j},
\label{eq:F25}
\end{equation}
where $\varphi_j$ are the eigenfunctions of $\cF_{\#}$ for the eigenvalues 
$\mu_j$. This function should be bounded if and only if $\vz \in \Omega$.

In practice, we can work only with the propagating part of the
scattered field, because the array is at large distance from the obstacle. Thus,  instead of $\cF$ defined as in Lemma
\ref{lem.1}, we use its projection on the subspace
\begin{equation}
{\mathscr P} = \mbox{span}\{\psi_0, \ldots, \psi_J\} \subset L^2(\cA).
\label{eq:F26}
\end{equation}
The projection is the $(J+1) \times (J+1)$ matrix 
\begin{equation}
\cF^{\csP} = \Big ( \big( \cF \psi_j, \psi_l \big)_{\cA} \Big)_{0 \le j,l \le J},
\label{eq:projF}
\end{equation}
which defines in turn the $(J+1) \times (J+1)$ Hermitian,
positive definite matrix 
\begin{equation}
\label{eq:propFsharp}
\cF^{\csP}_{\#} = \big |\Re (\cF^{\csP}) \big|  + { |\Im (\cF^{\csP})|}.
\end{equation}  
The implementation of the factorization method in section \ref{sect:numerics} is
based on the Picard range criterium for the square root of \eqref{eq:propFsharp}, so the series in \eqref{eq:F25}
becomes a finite sum with $J+1$ terms. The resulting image is expected to be larger
outside the obstacle, and the numerical results illustrate that this is indeed the case. However,  the equivalent of Theorem \ref{thm.Fsharp} is not yet established for the projection $\cF^{\csP}_{\#}$ to the propagating modes.

\section{Connection to migration imaging}
\label{sect:migration}
We describe in section \ref{sect:mig1}  the classic migration imaging function, where the scattered 
wave $\us$ is backpropagated to the search point $\vz$ using the Green's function in the empty waveguide. 
Then, we give in section \ref{sect:mig2} a slight modification of the migration imaging function, where the backpropagation 
is  done with the second 
derivative of the Green's function, for improved focusing of the image. The connection to  the factorization method is in section 
\ref{sect:mig3}.

\subsection{Migration imaging}
\label{sect:mig1}
Let $\PP:L^2(\cA) \to {\mathscr P}$ be the orthogonal
projector from $L^2(\cA)$ to $\csP$ and denote by
\begin{equation}
G_{\csP}(\cdot, \vz)\big|_{\cA} = \PP G(\cdot, \vz)\big|_{\cA}
\label{eq:R7}
\end{equation}
the propagating part of the Green's function evaluated at the array. The classic migration imaging function is given by 
\begin{equation}
{\cal J}(\vz) = \iint_{\cA} dS_{\vx} dS_{\vy} \, \us(\vx,\vy)  \overline{G_{\csP}(\vx,\vz)} \overline{G_{\csP}(\vy,\vz)},
\label{eq:Mi1}
\end{equation}
Because the array is far from the obstacles, we neglect the evanescent 
part of the measured $\us$ and backropagate it to $\vz$ using \eqref{eq:R7}.

Note from \eqref{eq:uinc} that $\overline{G_{\csP}(\cdot,\vz)\big|_{\cA}}$ is 
of the form \eqref{eq:F15}, so we can use \eqref{eq:F18}, the factorization \eqref{eq:F13} and the duality relation 
\eqref{eq:dualT} to rewrite \eqref{eq:Mi1} as 
\begin{equation}
{\cal J}(\vz) = \big(\cF \overline{G_{\csP}(\cdot,\vz)\big|_{\cA}}, \overline{G_{\csP}(\cdot,\vz)\big|_{\cA}}\big)_{\cA} = 
\left< \Lambda \cT  \overline{G_{\csP}(\cdot,\vz)\big|_{\cA}}, \cT  \overline{G_{\csP}(\cdot,\vz)\big|_{\cA}}\right>_{\partial \Omega}.
\label{eq:Mi2}
\end{equation}
We also obtain from definition \eqref{eq:F2} and the orthogonality relation \eqref{eq:EigFOrt}  that 
\begin{align}
K_0(\vx,\vz) &= \cT  \overline{G_{\csP}(\cdot,\vz)\big|_{\cA}}(\vx) =  \int_{\cA} dS_{\vy} \, G(\vx,\vy) \overline{G_{\csP}(\vy,\vz)} \nonumber \\
&=
\sum_{j=0}^J \frac{1}{\beta_j^2} \psi_j(\bxp) \psi_j(\bzp) 
\cos(\beta_j x) \cos(\beta_j z),  \qquad \vx = (x,\bxp) \in \partial \Omega.
\label{eq:Mi3}
\end{align}
In \eqref{eq:Mi2} we calculate the duality pairing 
\begin{equation}
\mathcal{J}(\vz) = \left< \Lambda K_0(\cdot,\vz) \big|_{\cA},  K_0(\cdot,\vz) \big|_{\cA}\right>_{\partial \Omega}
=  \int_{\partial \Omega}dS_{\vx} \,  \overline{h(\vx)} K_0(\vx,\vz),
\label{eq:Mi4}
\end{equation}
where $h = \Lambda K_0(\cdot,\vz) \big|_{\cA}$ is the solution of 
\begin{equation}
\int_{\partial \Omega} dS_{\vy} \, \overline{G(\vx,\vy)} h(\vy) = - K_0(\vx,\vz), \qquad 
\vx \in \partial \Omega.
\label{eq:Mi5}
\end{equation}
Because $\Lambda$ is an isomorphism,  we have that $\|h\|_{H^{-\frac{1}{2}}(\partial \Omega)}$  is large when $\|K_0(\cdot, \vz)\|_{H^\frac{1}{2}(\partial \Omega)}$ is large, so the focusing of the imaging function \eqref{eq:Mi4} depends on how sharply peaked the 
kernel $\eqref{eq:Mi3}$ is at $\vx = \vz$.

\begin{figure}
\centering
\includegraphics[width=0.52\textwidth]{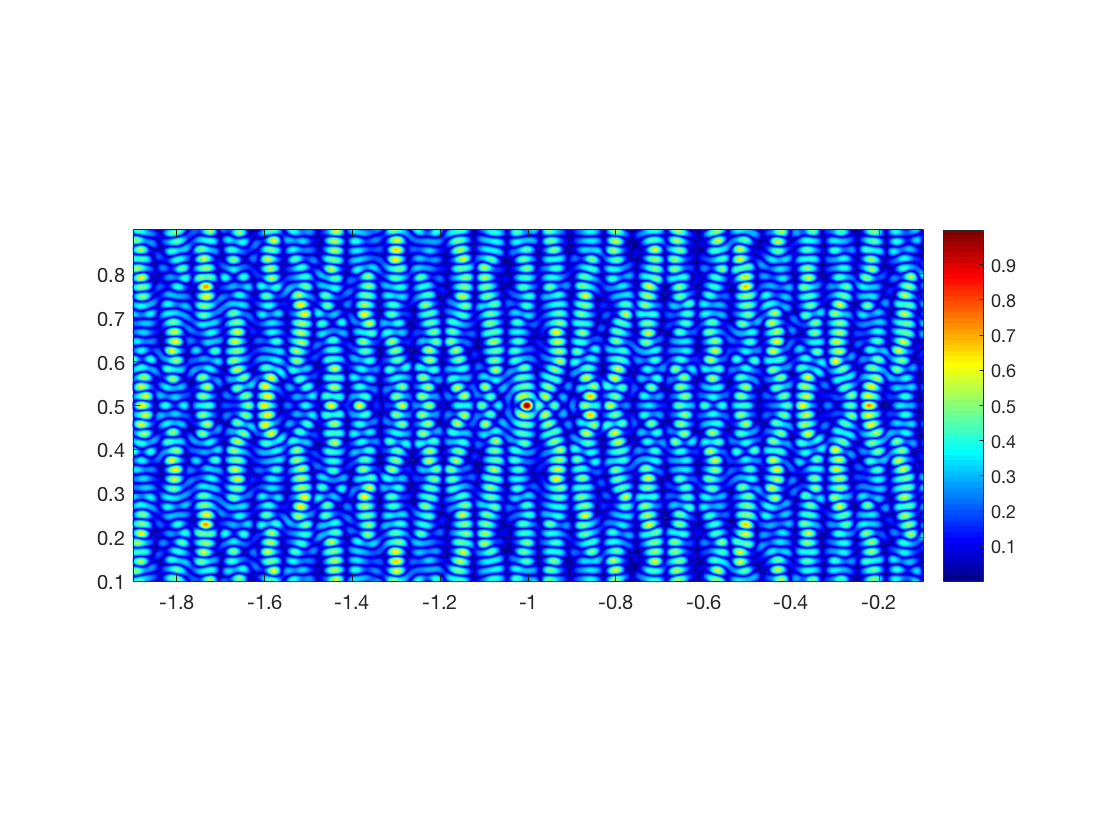}
\hspace{-0.3in}\includegraphics[width=0.52\textwidth]{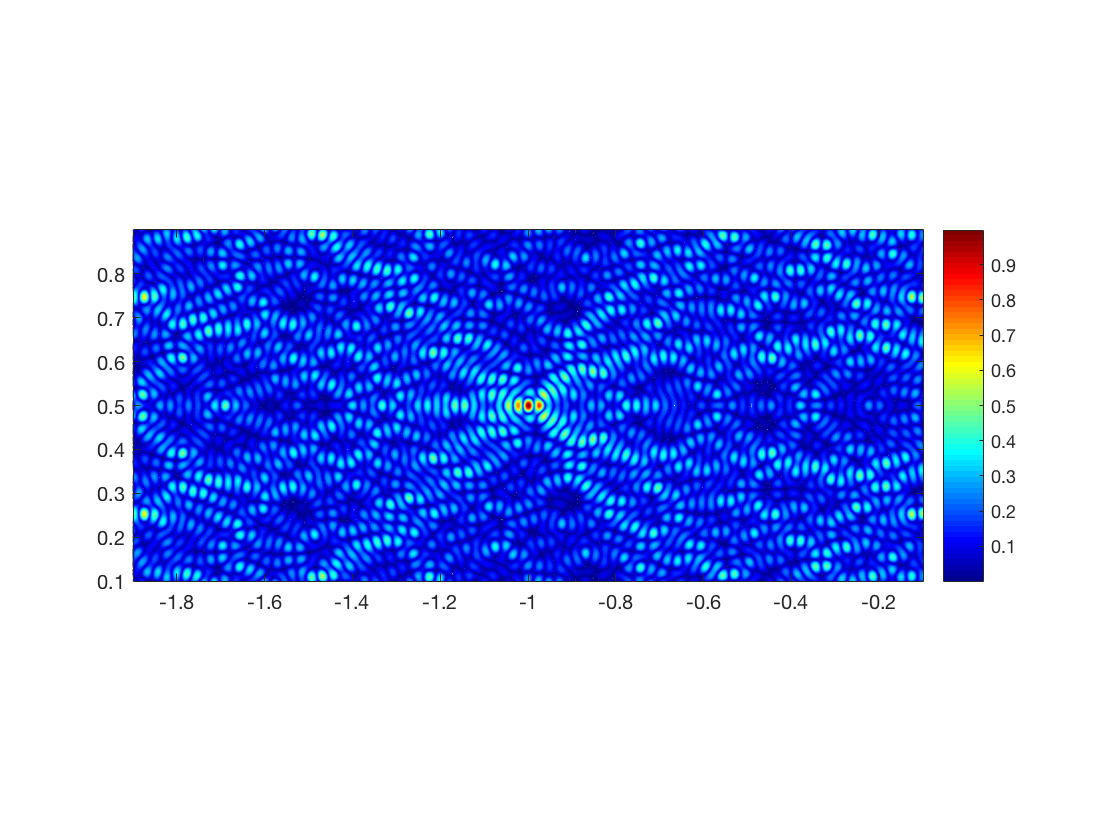}
\vspace{-0.6in}
\caption{The absolute value of the kernel $\cK_0(\vx,\vz)$ (left) and $\cK(\vx,\vz)$ (right)  in a two dimensional waveguide with
  $50$ propagating modes.  Both kernels are normalized by their maximum value. The point $\vx=(-|\cX|,|\cX|/2) $ is fixed
  and the search domain of $\vz =(z,\bzp)$ is indicated in the labels,
  in units of $|\cX|$. The abscissa is $z$ and ordinate is $\bzp$, in units of $|\cX|$.}
\label{fig:Kernels}
\end{figure}

We display in the left plot of Figure \ref{fig:Kernels} the kernel $K_0(\vx,\vz)$ in a two dimensional waveguide with $50$ propagating modes (see also Figure \ref{fig:KernelsH1}). We note that while $K_0(\vx,\vz)$  has 
a peak at $\vx = \vz$, there are many other peaks. In the next section we modify slightly the imaging function,
by backpropagating with the second range derivative  of $G_{\csP}$. This results in the better 
focused kernel $K(\vx,\vz)$ displayed in the right plot of Figure \ref{fig:Kernels}.

\subsection{A modified migration imaging function}
\label{sect:mig2}
Instead of using  $\overline{G_{\csP}(\cdot,\vz)\big|_{\cA}}$ to backpropagate the measured $\us$  to the imaging point $\vz$,  consider \begin{align}
\varphi_{\vz}(\vx) &= C_{\vz} \sum_{j=0}^J \frac{\beta_j}{i} 
\psi_j(\bzp)\psi_j(\bxp) e^{i \beta_j x_{\cA}} \cos (\beta_j z) \nonumber \\
&= - C_{\vz}  \partial_x^{2}
\overline{G_{\csP}(\vx,\vz)}\big|_{\vx \in \cA},  \qquad 
\vx = (x_{\cA},\bxp) \in \cA,
\label{eq:MS1}
\end{align}
where $C_{\vz}$ is a positive normalization constant so that 
\begin{equation}
\big(\overline{G(\cdot, \vz)},\varphi_{\vz} \big)_{\cA} = C_{\vz} \sum_{j=0}^J \psi_j^2(\bzp) \cos^2(\beta_j) = 1.
\label{eq:MS0}
\end{equation}
This function $\varphi_{\vz}$ is of the form \eqref{eq:F16}, so  we can calculate $\mathbb{P} \cF \varphi_{\vz}$ from the measurements at the array, using Lemma \ref{lem.1} and the matrix \eqref{eq:projF}.

\begin{figure}
\centering
\includegraphics[width=0.49\textwidth]{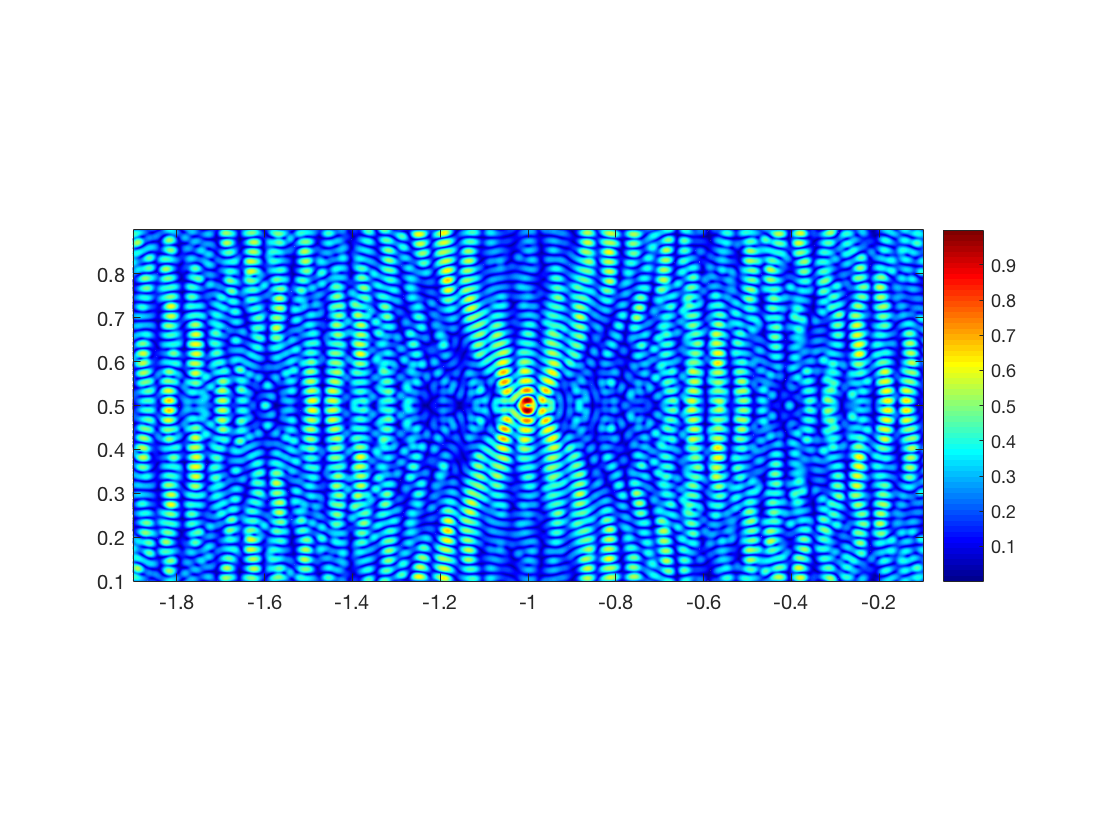}
\hspace{-0.3in}\includegraphics[width=0.49\textwidth]{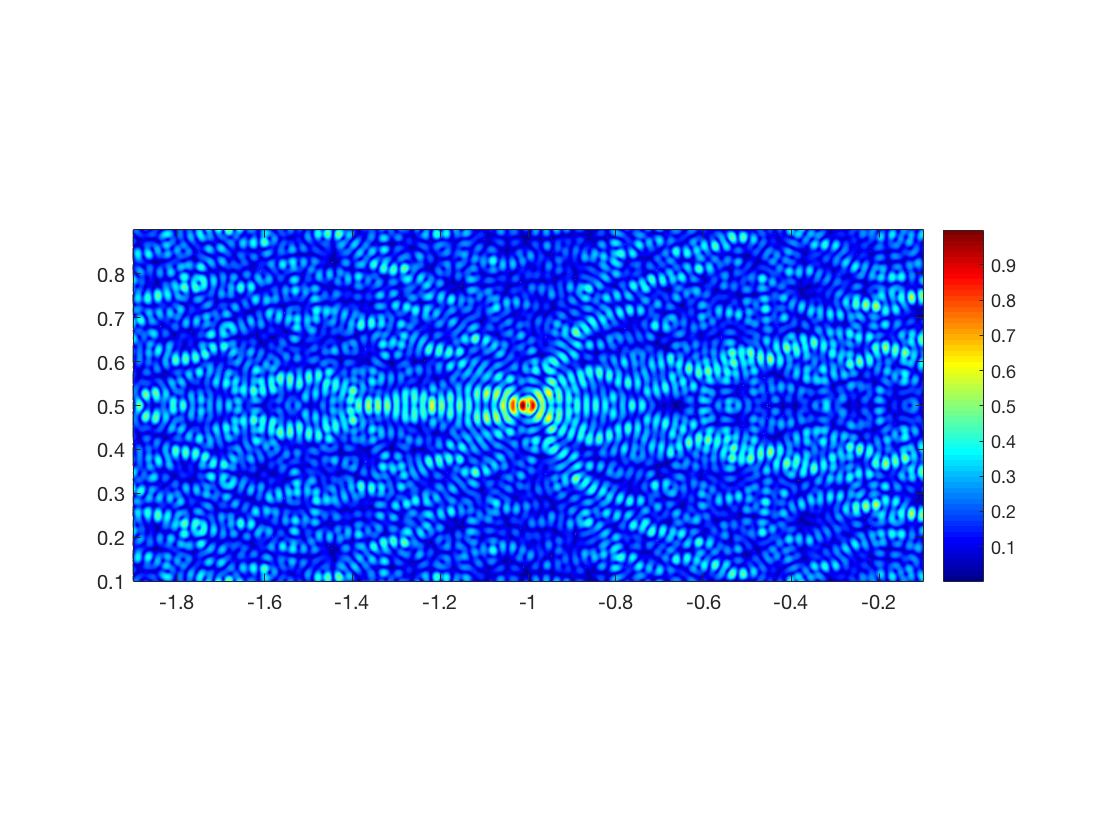}
\vspace{-0.6in}
\caption{We display  $\|\cK_0(\cdot,\vz)\|_{H^1(\Omega)}$ (left) and $\|\cK(\cdot,\vz)\|_{H^1(\Omega)}$ (right)  in a two dimensional waveguide with  $50$ propagating modes, where $\Omega$ is a square centered at $(-|\cX|,|\cX|/2)$, of side length $0.02|\cX|$.  In both plots we normalize to maximum value $1$.  The search domain of $\vz =(z,\bzp)$ is indicated in the labels,
  in units of $|\cX|$. The abscissa is $z$ and ordinate is $\bzp$, in units of $|\cX|$.}
\label{fig:KernelsH1}
\end{figure}

The modified migration type imaging function is 
\begin{equation}
\mathcal{J}_{\rm mig}(\vz) =  - \Im \Big[ \big(\cF \varphi_{\vz},\varphi_{\vz}
\big)_{\cA} \Big]= - \Im \Big[ \big(\mathbb{P}\cF \varphi_{\vz},\varphi_{\vz}
\big)_{\cA} \Big] = -\big(\Im(\cF) \varphi_{\vz},\varphi_{\vz}
\big)_{\cA},
\label{eq:Mi10}
\end{equation}
where we used the orthogonality relation \eqref{eq:EigFOrt}, definition \eqref{eq:F13p}  and the identity 
\[
\big(\cF^\star \varphi_{\vz},\varphi_{\vz}\big)_{\cA} = \big(\varphi_{\vz}, \cF \varphi_{\vz}\big)_{\cA}
= \overline{ \big(\cF \varphi_{\vz},\varphi_{\vz}
\big)_{\cA} }.
\]
We take the imaginary part in order to relate \eqref{eq:Mi10}  to the factorization method. 
Using  equation \eqref{eq:F13} in \eqref{eq:Mi10} we obtain 
\begin{equation}
\mathcal{J}_{\rm mig}(\vz) =  - \left< \Im (\Lambda) \cK(\cdot, \vz), \cK(\cdot, \vz) \right>_{\partial \Omega},
\label{eq:Mi11}
\end{equation}
where we introduced the kernel 
\begin{align}
K(\vx,\vz) &= \cT  \varphi_{\vz} (\vx) =  \int_{\cA} dS_{\vy} \, G(\vx,\vy) \varphi_{\vz}(\vy)
= C_{\vz}  \sum_{j=0}^J
\psi_j(\bxp)\psi_j(\bzp) \cos(\beta_j x) \cos( {\beta_j} z) \nonumber \\ &= C_{\vz}  \Re \Big[\partial_x
G_{\csP}(\vx,\vz)\Big], 
 \qquad \vx = (x,\bxp) \in \partial \Omega.
\label{eq:Mi12}
\end{align}
This kernel is peaked at $\vx = \vz$ and decays with $|\vx-\vz|$ as illustrated in the right plots of Figures 
\ref{fig:Kernels} and \ref{fig:KernelsH1}.  

Because $\Lambda \cK(\cdot, \vz)\big|_{\partial \Omega}$ is bounded in $H^{\frac{1}{2}}(\partial \Omega)$, 
the imaging function \eqref{eq:Mi11} is bounded above in terms of $ \|\cK(\cdot, \vz)\|_{H^{\frac{1}{2}}(\partial \Omega)}$ and therefore of 
$ \|\cK(\cdot, \vz)\|_{H^{1}( \Omega)}$. The latter norm is small when $\vz$ is far from $\Omega$, as illustrated in 
Figure \ref{fig:KernelsH1}. 
By Theorem \ref{thm.1}, the operator $-\Im(\Lambda)$ is self-adjoint and positive semi-definite, so we 
expect that the imaging function \eqref{eq:Mi11} is large for search points $\vz$ near $\partial \Omega$, as long as 
$K(\cdot,\vz)$ is not in the null space of $\Im(\Lambda)$. 

The next theorem sheds more light on the  behavior of $\mathcal{J}_{\rm mig}(\vz)$ for search points near 
$\Omega$. To state it, let  $\itbf{U}^{{\rm sc},\csP}$ be  the $(J+1) \times (J+1)$
matrix obtained by projecting the measured scattered field $\us$ on the
finite dimensional subspace \eqref{eq:F26}. The entries of this matrix
are
\begin{equation}
\UsP_{j,l} = \iint_{\cA} dS_{\vx} dS_{\vy} \, \us(\vx,\vy)
\psi_j(\bxp) \psi_l(\byp), \qquad j,l = 0, \ldots, J,
\label{eq:R1}
\end{equation}
and we note that $\itbf{U}^{{\rm sc},\csP}$ is complex symmetric, by
reciprocity, but it is not Hermitian. The singular value decomposition
of $\itbf{U}^{{\rm sc},\csP}$ is of the form
\begin{equation}
\itbf{U}^{{\rm sc},\csP} = \itbf{V} \boldsymbol{\mathfrak{S}} \itbf{V}^T,
\label{eq:R2}
\end{equation}
where ${\itbf V}$ is unitary, with columns $\bv_j$ for $j = 1, \ldots,
J+1$, and $\boldsymbol{\mathfrak{S}}$ is the diagonal matrix of
singular values, in decreasing order. Typically, the matrix $\itbf{U}^{{\rm sc},\csP}$ 
is rank defficient, with rank  $r < J+1$. Its null space is spanned by
the right singular vectors $\overline{\bv_j}$, for $j = r+1, \ldots,
J+1$. We denote by $\overline{V_{j,l}}$ the entries of these singular
vectors, and use them to define the following subspace of $\csP$, of
dimension $J-r+1$,
\begin{equation}
\csP_0 = \mbox{span} \Big\{ \sum_{j=0}^J
\overline{V_{j+1,l}}\psi_j(\bxp), ~~ l = r+1, \ldots, J+1\Big\}.
\label{eq:R4}
\end{equation}
This is in the null space of the operator $\PP \cN:L^2(\cA) \to
{\mathscr P}$. The orthogonal complement of
$\csP_0$ in $\csP$ is denoted by $\csP_0^\perp$, so we can write
\begin{equation}
\csP = \csP_0 \oplus \csP_0^\perp.
\label{eq:R5}
\end{equation}
The following theorem is proved in Appendix \ref{ap:D}.

\vspace{0.05in}
\begin{theorem}
\label{thm.NF}
Consider a search point $\vz \in (x_{\cA},0) \times \cX$, so that $\overline{G_{\csP}(\cdot, \vz)\big|_{\cA}} \notin \csP_0$. If $\vz \in \Omega$, then 
\begin{equation}
\mbox{inf}\Big\{ \big( { -\Im(\cF)} \varphi,\varphi \big)_{\cA}:
~ \varphi \in \csP_0^\perp, ~ ~ \big(\overline{G(\cdot,
  \vz)},\varphi \big)_{\cA} = 1 \Big\} > 0.
\label{eq:R6}
\end{equation}
\end{theorem}

\vspace{0.05in} This result,  the factorization \eqref{eq:F13} and definition \eqref{eq:F13p} imply  that 
when $\vz \in \Omega$, we have 
\[
 \big( { -\Im(\cF)} \varphi,\varphi \big)_{\cA} = - \left< \Im (\Lambda) \cT \varphi, \cT \varphi \right>_{\partial \Omega} > 0,
 \]
for all $\varphi \in \csP_0^\perp$ normalized by
$
\big(\overline{G(\cdot, \vz)},\varphi \big)_{\cA}  = 1.
$
The function $\varphi_{\vz}$ defined in \eqref{eq:Mi10} satisfies this normalization but it may not lie in $\csP_0^\perp$. 
Thus,  there can be  points $\vz \in \Omega$ where $\mathcal{J}_{\rm mig}(\vz)$ is small. 

Theorem \ref{thm.NF} suggests another modification of  the migration imaging function, where 
the backpropagation is carried out with the projection of $\varphi_{\vz}$ on $\csP_0^\perp$. We do not consider such a  
modification in this paper, but introduce instead a new imaging function that is guaranteed not to vanish at $\vz \in \Omega$ and is  related to the formulation  \eqref{eq:FR6} of the factorization method.  

\subsection{Connection to the factorization method}
\label{sect:mig3}
The new migration type  imaging function backpropagates with   the same $\varphi_{\vz} \in \csP$ defined in \eqref{eq:MS1},
\begin{equation}
\mathcal{J}_{{\rm mig}\#}(\vz) =   \big(\cF_{\sharp} \varphi_{\vz},\varphi_{\vz}
\big)_{\cA} =  \big(\mathbb{P}\cF_{\sharp} \varphi_{\vz},\varphi_{\vz}
\big)_{\cA},
\label{eq:Mi20}
\end{equation}
where the last equality is due to the orthogonality relation \eqref{eq:EigFOrt}. It
can be computed from the array measurements using the matrix \eqref{eq:propFsharp}, and we can rewrite it 
using the factorization \eqref{eq:F23-1} and equation \eqref{eq:Mi12}, 
\begin{align}
\mathcal{J}_{{\rm mig}\#}(\vz) = \left< \Lambda_{\sharp} \cT \varphi_{\vz}, \cT \varphi_{\vz}\right>_{\partial \Omega} = 
\left< \Lambda_{\sharp} K(\cdot, \vz), K(\cdot, \vz)\right>_{\partial \Omega} \ge \mathcal{J}_{{\rm mig}}(\vz),
 \label{eq:Mi21}
\end{align}
where the inequality follows from equations  \eqref{eq:F13pp} and \eqref{eq:Mi11}.

The advantage of this imaging function is that the operator $\Lambda_{\sharp}$ is positive definite. As in the previous section, we expect that 
$\mathcal{J}_{{\rm mig}\#}(\vz)$ is large near the obstacle, due to the focusing property of the kernel $K(\vx,\vz)$, for $\vx \in \partial \Omega$. In fact, if $\mathcal{J}_{{\rm mig}}(\vz)$ is large at a point $\vz$, then $\mathcal{J}_{{\rm mig}\#}(\vz) $ is even larger. In addition, we can use Theorem \ref{thm.Fsharp} to conclude that since $\varphi_{\vz}$ is in the admissible set of the optimization in \eqref{eq:FR6}, we have 
\begin{equation}
\mathcal{J}_{{\rm mig}\#}(\vz) > 0, \qquad \forall \vz \in \Omega.
\label{eq:Mi22}
\end{equation}
For points $\vz \notin \Omega$, the imaging function decays with the distance from $\vz$ to $\partial \Omega$,  because of the decay of $\|K(\cdot,\vz)\|_{H^1(\Omega)}$ illustrated in Figure \ref{fig:KernelsH1}.

Note that in theory,  the factorization method 
should perform better than the migration type imaging function, because in Theorem \ref{thm.Fsharp} we minimize  
$\big(\cF_{\sharp} \varphi,\varphi \big)_{\cA}$ over all the test functions $\varphi$ in \eqref{eq:R6}, whereas in 
\eqref{eq:Mi20} we consider a single test function $\varphi_{\vz}$.  However, the migration method has
the advantage that it combines easily multiple frequency measurements,
by simply superposing \eqref{eq:Mi20} at the given
frequencies. This results in a significant improvement of the 
images, as illustrated in section \ref{sect:numerics}. To our
knowledge, there is no satisfactory way to take advantage of multiple
frequency data in the factorization method. The numerical results in section \ref{sect:numerics} also illustrate that the migration imaging function is more robust to noise and limited array aperture.

\section{Numerical results}
\label{sect:numerics}
In this section we present a comparative numerical study of the factorization and migration 
imaging methods in two dimensions.  

In the simulations, all lengths are in units of  $|\cX|$, the length 
of the cross-section interval $\cX = (0,|\cX|)$.  The scattered field  $\us$ is obtained by solving the wave equation in the sector 
$(-5|\cX|,0) \times \cX$ of the waveguide, using the high-performance multi-physics finite element 
software Netgen/NGSolve \cite{schoberl1997netgen} and a perfectly matched layer at range $-5|\cX|$. The array response matrix $ \itbf{U}^{\rm sc}$ defined in \eqref{eq:RESPM} is obtained by sampling 
$\us(\vx_r,\vx_s)$ at equidistant points  in $\cA = \{-2|\cX|\} \times \cX$, separated by $|\cX|/60$. 
It is contaminated with additive, complex Gaussian, iid noise with {standard deviation $\sigma_{\rm noise}$  calculated as a percent of the maximum 
absolute value of the entries in $\itbf{U}^{\rm sc}$}. 

We work only with the propagating modes, so we transform $\itbf{U}^{\rm sc}$ to the matrix 
$\itbf{U}^{{\rm sc},\csP} \in \mathbb{C}^{(J+1) \times (J+1)}$ defined in \eqref{eq:R1}, using the eigenfunctions  
\begin{equation}
\psi_j(\bxp) = \sqrt{\frac{2-\delta_{j,0}}{|\cX|}} \cos \Big(\frac{j \pi \bxp}{|\cX|}\Big).
\end{equation}
The integrals in \eqref{eq:R1} are approximated by Riemann sums, using the discrete sample points in $\cA$. 

We present results for two wavenumbers: $k = 29.15 \pi/|\cX|$ and $k = 49.15 \pi/|\cX|$, so that the 
waveguide supports $J+1 = 30$ and $50$ propagating modes, respectively.  For the migration images we also present multifrequency
results obtained at the wavenumbers $(29 + 0.15 m)\pi/|\cX|$,  with $m = 1, \ldots, 6$. The imaging region swept 
by the search point $\vz$ is $(-1.9 |\cX|, -0.1|\cX|) \times (0.1|\cX|, 0.9|\cX|)$.

To assess how the size of the array aperture affects the quality of the images, we present full and partial aperture results, 
where the array lies in the set $\{-2|\cX|\} \times (0, |\cX|_{\cA})$, with $|\cX|_{\cA} \le |\cX|$.
The implementation of the migration method is independent of the size of the aperture. For the factorization method and the modified 
migration method \eqref{eq:Mi12} 
we first process the partial aperture data as explained in \cite[Section 2.4]{borcea2018direct}, in order to obtain an estimate of 
the  matrix $\itbf{U}^{{\rm sc},\csP}$ used in Algorithms \ref{alg.1}--\ref{alg.2} below. The migration method \eqref{eq:Mi10} 
calculated in Algorithm \ref{alg.3} does not require this extra data processing.

\subsection{Imaging algorithms}

The implementation of the factorization method is as  described in 
section \ref{sect:fact.2}, except that we use only  the propagating part of the data:

\vspace{0.05in}
\begin{algorithm} The factorization method:
\label{alg.1}

\vspace{0.05in}
\noindent \textbf{Input:} The  matrix $\itbf{U}^{{\rm sc},\csP}$ (with or without noise) and the imaging mesh.

\noindent \textbf{Processing steps:}
\vspace{0.05in}
\begin{enumerate}
\item Represent the operator $\cF$ by the $(J+1) \times (J+1)$ matrix 
\[
\cF^{\csP} = \Big ( \big( \cF \psi_j, \psi_l \big)_{\cA} \Big)_{j,l = 0, \ldots, J} = \Big(  -\overline{\itbf{U}_{jl}^{{\rm sc},\csP}} e^{-i 2\beta_l x_\cA} \Big)_{j,l = 0, \ldots, J},
\]
where we used Lemma \ref{lem.1} and equation \eqref{eq:R1}.
\item Calculate the matrix $\cF^{\csP}_{\#} = \big |\Re (\cF^{\csP}) \big|  + { |\Im (\cF^{\csP})|}$, which is Hermitian, 
positive definite,  with the eigenvalue decomposition  $\cF^{\csP}_\# = {\bf V} {\bf D} {\bf V}^\star,$
where the star denotes complex conjugate and transpose.  Its square root is 
$ ( \cF^{\csP}_\# )^{\frac{1}{2}} = {\bf V} {\bf D}^{\frac{1}{2}} {\bf V}^\star.$
Denote by $\bv_j$ the columns of the unitary matrix ${\bf V}$ and by $d_{jj} \ge 0$ the entries of 
${\bf D}$, for $j = 0, \ldots, J+1$.
\item For all $\vz $ on the imaging mesh and a user defined small parameter $\epsilon >0$ calculate 
the regularized solution ${\bf g}^\epsilon_{\vz}$ of $( \cF^{\csP}_\# )^{\frac{1}{2}} {\bf g}_{\vz}  = {\bf b}_{\vz}$,
where ${\bf b}_{\vz} \in \CC^{J+1}$ is the  column vector with entries
\[
b_{j,\vz} = \int_{\cA} dS_{\vx} \, \overline{G_{\csP}(\vx,\vz)} \psi_j(\bxp),
\qquad j = 0, \ldots, J.
\]
This regularized solution satisfies 
\begin{eqnarray*}
\|{\bf g}^\epsilon_{\vz} \|^2 = \sum_{j=0}^{J} |{\bf b}^\star_{\vz} \bv_j|^2\frac{d_{jj}}{(d_{jj} + \alpha^\epsilon)^2},
\end{eqnarray*}
where  $\alpha^\epsilon$ is a positive Tikhonov regularization parameter chosen according to the Morozov principle, so that
\begin{eqnarray*}
\| ( \cF^{\csP}_\# )^{\frac{1}{2}}  {\bf g}^\epsilon_{\vz}  - {\bf b}_{\vz}  \| = \epsilon \| {\bf g}^\epsilon_{\vz}  \|.
\end{eqnarray*}
\item Calculate the imaging function
\[
{\cal J}_{\#}(\vz) =  \frac{1/\|{\bf g}^\epsilon_{\vz} \|}{\sup_{\vz'} 1/\|{\bf g}^\epsilon_{\vz'} \|}.
\]
\end{enumerate}

\noindent \textbf{Output:} The estimate of the support of $\Omega$ is
determined by the set of points $\vz$ where ${\cal J}_\#(\vz)$ is larger than the user defined threshold. 
\end{algorithm}

\vspace{0.05in}
The   migration type imaging function is  \eqref{eq:Mi12} calculated with the following algorithm:
\vspace{0.05in}
\begin{algorithm}  Imaging with $\mathcal{J}_{{\rm mig}\#}(\vz)$:
\label{alg.2}

\noindent \textbf{Input:} The matrix $\itbf{U}^{{\rm sc},\csP}$  (with or without noise) and the imaging mesh.

\noindent \textbf{Processing steps:}
\vspace{0.05in}
\begin{enumerate}
\item Calculate $\cF^{\csP}$ and $\cF_{\sharp}^{\csP}$ as in Algorithm \ref{alg.1}. 
\item For all $\vz$ on the imaging mesh, calculate the  column vector $\boldsymbol{a}_{\vz} \in \CC^{J+1}$, 
with entries 
\[
a_{j,\vz} = \int_{\cA} dS_{\vx} \, \varphi_{\vz}(\vx_r) \psi_j(\bxp),
\qquad j = 0, \ldots, J,
\]
where $\varphi_{\vz}$ is defined in  \eqref{eq:MS1}.
\item 
Calculate 
\begin{align*}
\mathcal{J}_{{\rm mig}\#}(\vz) &= \boldsymbol{a}_\bz^\star \cF_{\sharp}^{\csP} \boldsymbol{a}_\bz 
\quad \mbox{and then} \quad \mathcal{J}_{{\rm mig}\#}(\vz) &= \frac{\mathcal{J}_{{\rm mig}\#}(\vz)}{\max_{\vz'} \mathcal{J}_{{\rm mig}\#}(\vz')},
\end{align*}
where the star denotes complex conjugate and transpose. 
\end{enumerate}

\noindent \textbf{Output:} The estimate of the support of $\Omega$ is
determined by the set of points $\vz$ where $\mathcal{J}_{{\rm mig}\#}(\vz) $  is  larger than the user defined threshold. 
\end{algorithm}

\vspace{0.05in}
The   migration imaging function  \eqref{eq:Mi10} is calculated with the following algorithm:
\vspace{0.05in}
\begin{algorithm}  Imaging with $\mathcal{J}_{{\rm mig}}(\vz)$:
\label{alg.3}

\noindent \textbf{Input:} The $n_{\cA} \times n_{\cA}$ array response matrix $\itbf{U}^{{\rm sc}}$ defined in \eqref{eq:RESPM} (with or without noise) and the imaging mesh.

\noindent \textbf{Processing steps:}
\vspace{0.05in}
\begin{enumerate}
\item For all $\vz$ on the imaging mesh, calculate the  column vector $\boldmath{\phi}_{\vz} \in \CC^{n_{\cA}}$, 
with entries defined by $\varphi_{\vz}$ evaluated at the sensor locations $\vx_r$, 
\[
\phi_{r,\vz} = \varphi_{\vz}(\vx_r), 
\qquad r = 1, \ldots, n_{\cA}.
\]
\item 
Calculate 
\begin{align*}
\mathcal{J}_{{\rm mig}}(\vz) &= \big| \Im \big(\boldmath{\phi}_{\vz}^T \itbf{U}^{{\rm sc}} \boldmath{\phi}_{\vz}\big) \big| 
\quad \mbox{and then} \quad 
\mathcal{J}_{{\rm mig}}(\vz) &= \frac{\mathcal{J}_{{\rm mig}}(\vz)}{\max_{\vz'} \mathcal{J}_{{\rm mig}}(\vz')}.
\end{align*}
\end{enumerate}

\noindent \textbf{Output:} The estimate of the support of $\Omega$ is
determined by the set of points $\vz$ where $\mathcal{J}_{{\rm mig}}(\vz) $  is  larger than the user defined threshold. 
\end{algorithm}

\subsection{Numerical results}
\begin{figure}[t]
\includegraphics[width=1\textwidth]{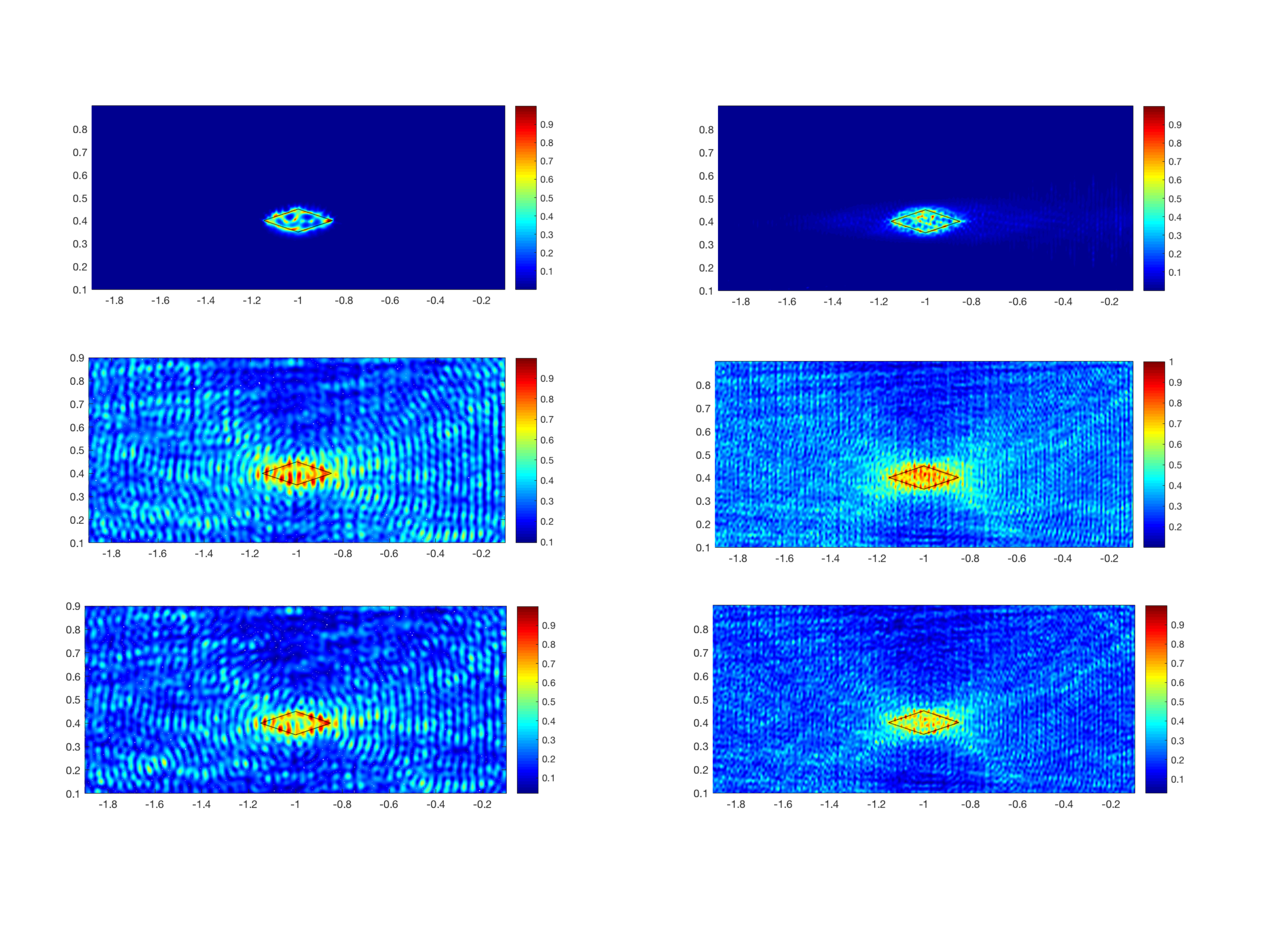} 
\vspace{-1.9cm}
\caption{Reconstruction of a rhombus shaped obstacle shown
  with a solid black line.  The abscissa is range and the ordinate is
  cross-range, scaled by $|\cX|$. Full aperture, noiseless array data. Top line: $\mathcal{J}_{\#}(\vz)$ . Middle line: 
  $\mathcal{J}_{{\rm mig}\#}(\vz) $. Bottom line: $\mathcal{J}_{{\rm mig}}(\vz)$. Left column: 30 propagating modes. Right column: 50 propagating modes.}
\label{fig:FactMig_FA}
\end{figure}

We now present results obtained with Algorithms \ref{alg.1}--\ref{alg.3}. In Figure \ref{fig:FactMig_FA} we display the effect
of the  probing frequency and therefore of the number of propagating modes. As expected, the higher the 
frequency, the better the resolution. The remaining images in this section are obtained in a waveguide with $30$ propagating modes.

\begin{figure}[h]
\centering 
\includegraphics[width=1\textwidth]{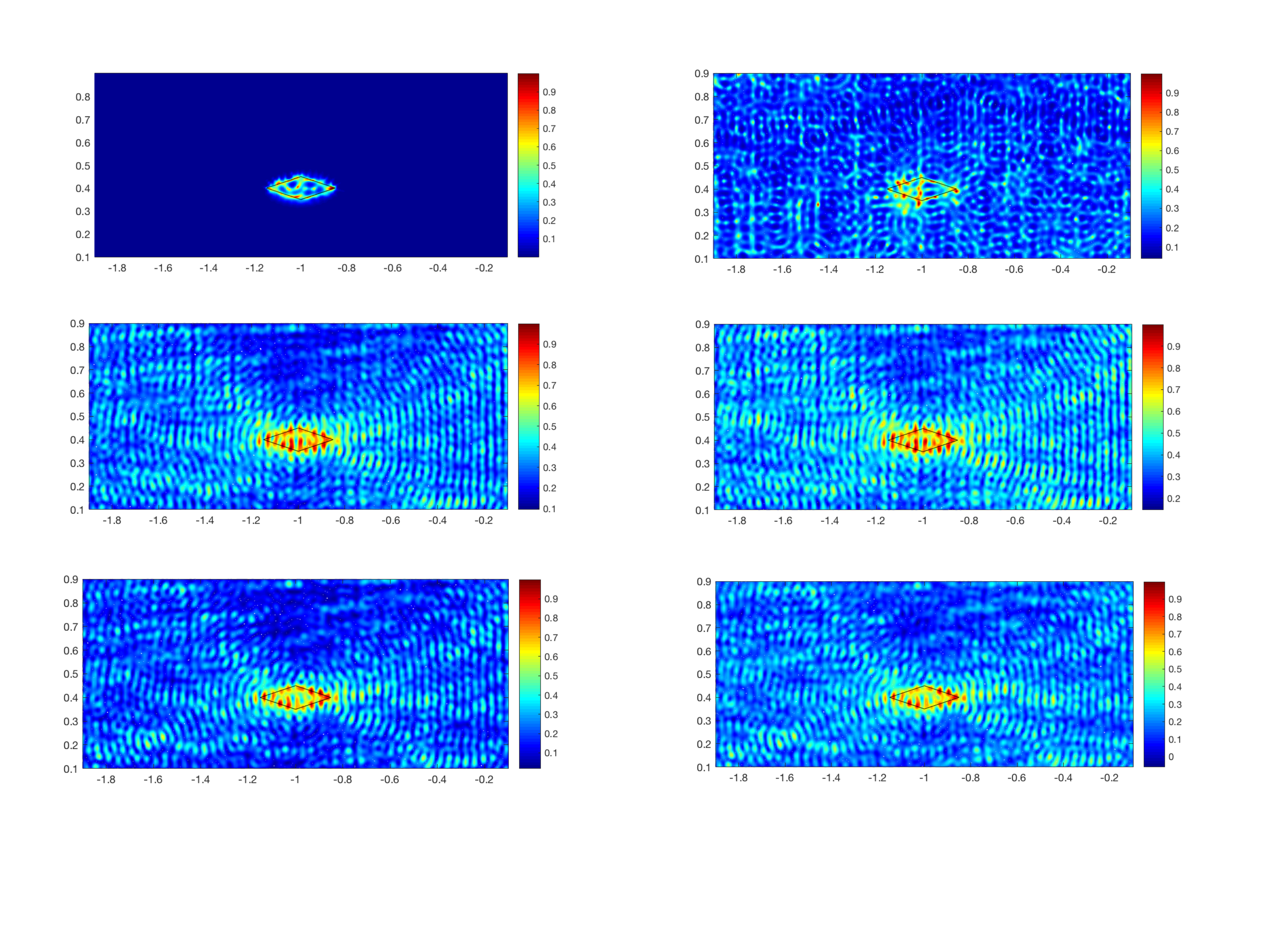}
\vspace{-2.2cm}\caption{Reconstruction of a rhombus shaped obstacle shown
  with a solid black line.  The abscissa is range and the ordinate is
  cross-range, scaled by $|\cX|$. Full aperture array data and $30$ propagating modes. Top line: $\mathcal{J}_{\#}(\vz)$ . Middle line: 
  $\mathcal{J}_{{\rm mig}\#}(\vz) $. Bottom line: $\mathcal{J}_{{\rm mig}}(\vz)$. Left column: no noise. Right column: $10\%$ noise. }
\label{fig:Noise comparison}
\end{figure}

\begin{figure}[h]
\centering 
\includegraphics[width=1\textwidth]{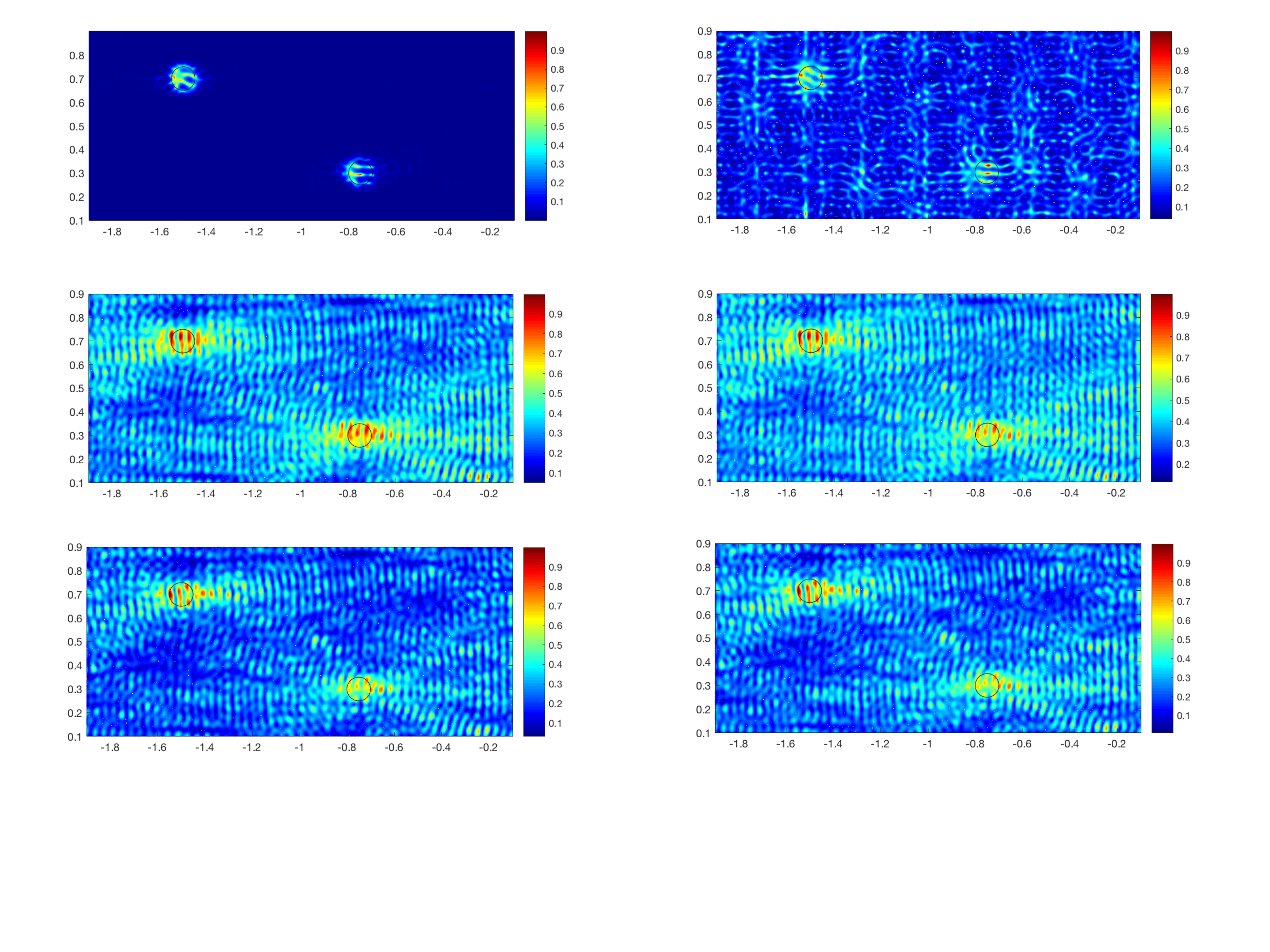}  
\vspace{-2.5cm}\caption{Reconstruction of two obstacles.  The abscissa is range and the ordinate is
  cross-range, scaled by $|\cX|$. Full aperture array data and $30$ propagating modes. Top line: $\mathcal{J}_{\#}(\vz)$ . Middle line: 
  $\mathcal{J}_{{\rm mig}\#}(\vz) $. Bottom line: $\mathcal{J}_{{\rm mig}}(\vz)$. Left column: no noise. Right column: $10\%$ noise.}
\label{fig:CIRCLE_FA}
\end{figure}

The robustness to noise is illustrated in Figures \ref{fig:Noise comparison} and \ref{fig:CIRCLE_FA}, where we display images 
of a rhombus shaped obstacle and two circle shaped obstacles obtained with noiseless data (left columns) and 
data contaminated with $\sigma_{\rm noise} = 10\%$ noise (right columns). 


\begin{figure}[h]
\begin{center}
\includegraphics[width=1\textwidth]{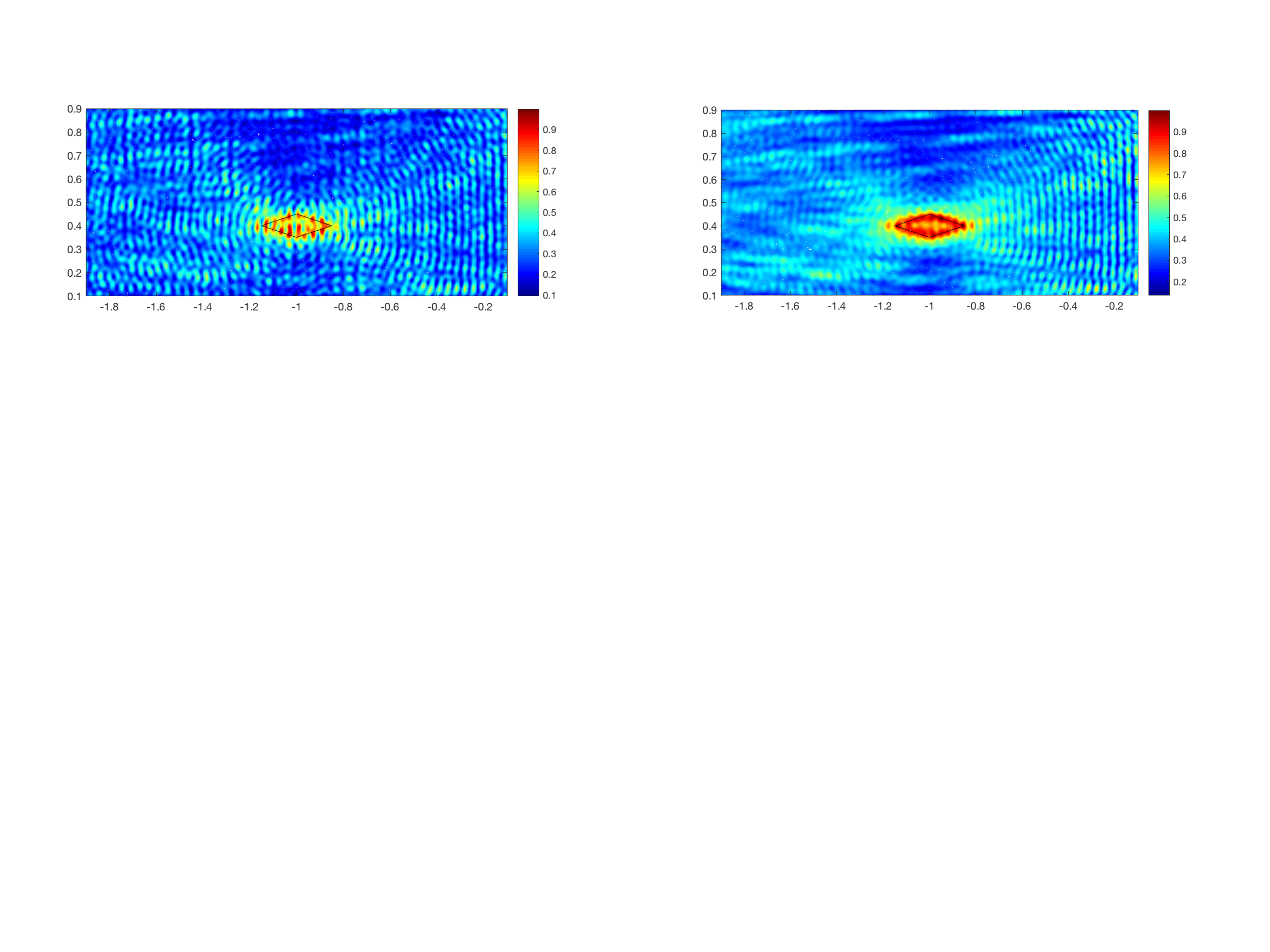}
\vspace{-7.0cm}
\end{center}
\begin{center}
\includegraphics[width=1\textwidth]{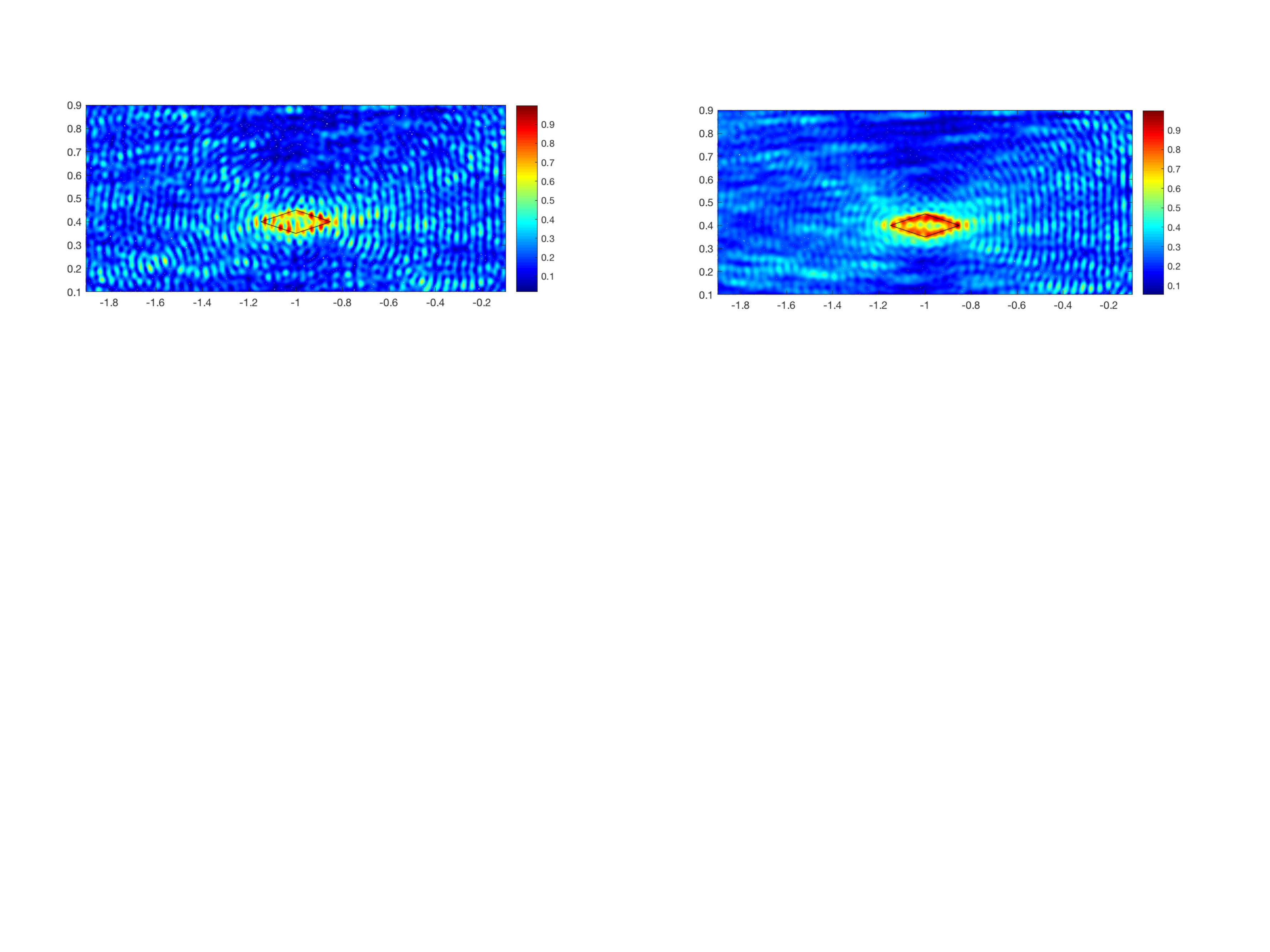}
\vspace{-7cm}\caption{Reconstruction of a rhombus shaped obstacle shown
  with a solid black line.  The abscissa is range and the ordinate is
  cross-range, scaled by $|\cX|$. Full aperture, noiseless array data and $30$ propagating modes. Top line: 
  $\mathcal{J}_{{\rm mig}\#}(\vz) $. Bottom line: $\mathcal{J}_{{\rm mig}}(\vz)$. Single frequency result (left) and multiple frequency result (right).}
\label{fig:Multi_Freq}
\end{center}
\end{figure}

In the noiseless case,  the results in Figures \ref{fig:FactMig_FA}--\ref{fig:CIRCLE_FA} show that the factorization 
method gives better images, as expected from the discussion at the end of section \ref{sect:mig3}. However, the 
migration images are most robust to noise i.e., they are similar for noiseless and the noisy data. Moreover, they improve 
significantly when we use multifrequency data, as illustrated in Figure  \ref{fig:Multi_Freq}.

\begin{figure}[h]
\centering 
\includegraphics[width=1\textwidth]{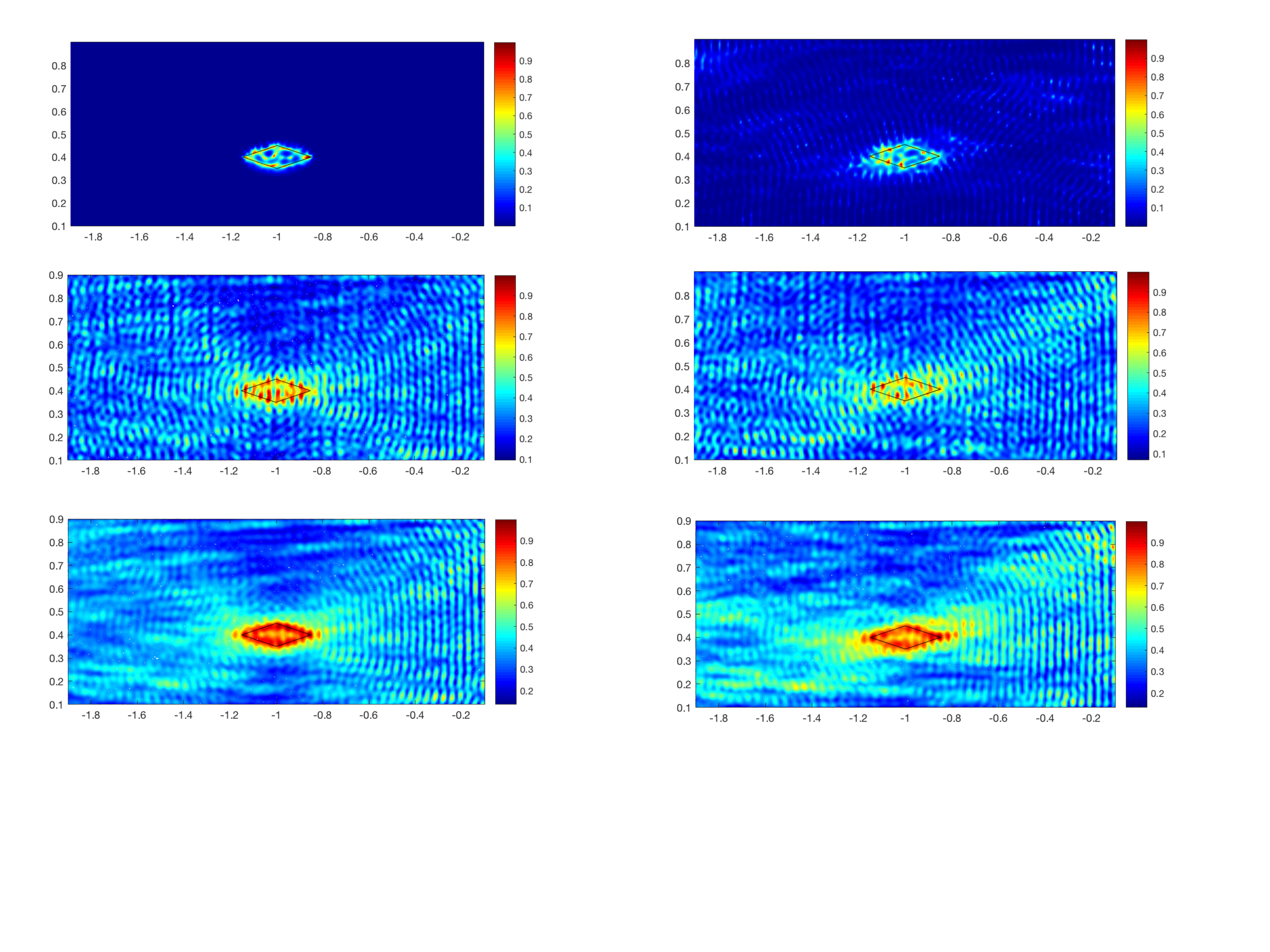} 
\vspace{-2.8cm}\caption{Reconstruction of a rhombus shaped obstacle shown
  with a solid black line.  The abscissa is range and the ordinate is
  cross-range, scaled by $|\cX|$. Noiseless array data. Top line: $\mathcal{J}_{\#}(\vz)$ . Middle  line: 
  $\mathcal{J}_{{\rm mig}\#}(\vz) $. Bottom line: $\mathcal{J}_{{\rm mig}_\#}(\vz)$ with multifrequency data. Left column: full aperture. 
  Right column: $75\%$ aperture array.}
\label{fig:partAp}
\end{figure}

The last images, in Figure \ref{fig:partAp} show the effect of the limited array aperture. They are obtained with $30$ propagating 
modes for noiseless data collected on an array of $|\cX|_{\cA} = 0.75 |\cX|$ aperture.  The images deteriorate at partial aperture,
but the migration method is clearly better when we use the multifrequency data.

\section{Summary} We presented a theoretical and computational comparative study of two qualitative methods for imaging 
obstacles in a terminating waveguide.  The first method is based on the factorization of the far field operator, defined by 
measurements of the scattered wave collected by an active array of sensors. It is designed to image at a single frequency and determines the support of the obstacles by either solving an optimization problem or, equivalently, using a Picard range criterium. 
The second method, known as migration, is based on the backpropagation of the  measured scattered wave to imaging points, using the Green's function in the empty waveguide. We studied the classic migration imaging method and explained how to modify it 
to get better images. Then, we related the migration type imaging method  to the factorization method, and  compared 
their performance with numerical simulations.

\label{sect:summary}
\section*{Acknowledgments}
This material is based upon research supported in part by the Air
Force Office of Scientific Research under award FA9550-18-1-0131. Part
of the research was done at the Institute for Computational and
Experimental Research in Mathematics in Providence, RI, during the
Fall 2017 semester. The authors' participation was supported by the
National Science Foundation under Grant No. DMS-1439786 and the Simons
Foundation Institute Grant Award ID 507536.

\appendix
\section{Proof of Lemma \ref{lem.1}}
\label{ap:A}
Consider first a function $\phi \in L^2(\cA)$ of the form
\begin{equation} 
\phi(\vx) = \sum_{j=0}^J \frac{\alpha_j}{i} \psi_j(\bxp) e^{i \beta_j
  x_{\cA}} + \sum_{j > J} \alpha_j \psi_j(\bxp), \qquad \vx \in \cA,
\label{eq:A1}
\end{equation}
with real valued coefficients $\alpha_j$, for $j \ge 0$, and  note from
the expression \eqref{eq:uinc} of the Green's function and the
orthogonality relation \eqref{eq:EigFOrt} that
\begin{align}
\int_{\cA} dS_{\vy} \,G(\vx,\vy) \phi(\vy) &= \sum_{j=0}^J
\frac{\alpha_j}{\beta_j} \psi_j(\bxp) \cos(\beta_j x) + \sum_{j > J}
\frac{\alpha_j}{|\beta_j|} \psi_j(\bxp) \cosh(|\beta_j| x) \nonumber
\\ &= \int_{\cA} dS_{\vy} \,\overline{G(\vx,\vy) \phi(\vy)}, \qquad
\vx \in \partial \Omega. \label{eq:A2}
\end{align}
Let us define,
\begin{align}
w(\vx) &= \int_{\cA} dS_{\vy} \, \us(\vx,\vy) \phi(\vy), \qquad \vx
\in \cW \setminus \overline{\Omega}, \label{eq:A3} \\ v(\vx) &= \int_{\cA}
dS_{\vy} \, \overline{\us(\vx,\vy) \phi(\vy)}, \qquad \vx \in \cW
\setminus \overline{\Omega}, \label{eq:A4}
\end{align}
and obtain from (\ref{eq:usc1}--\ref{eq:usc3}) and definitions 
(\ref{eq:F6}--\ref{eq:F8}) that
\begin{align}
v,w &\in \mathscr{H}(W\setminus \overline{\Omega}), \label{eq:A5} \\
w|_{\cA} &\in \mathscr{H}^{\rm out}(\cA), \label{eq:A6} \\
v|_{\cA} &\in \mathscr{H}^{\rm in}(\cA). \label{eq:A7}
\end{align}
Since $\us(\cdot,\vx_s)|_{\partial \Omega} = -
G(\cdot,\vx_s)|_{\partial \Omega}$, we also conclude from
(\ref{eq:A2}--\ref{eq:A4}) that
\begin{equation}
w(\vx) = v(\vx), \qquad \vx \in \partial \Omega.
\label{eq:A8}
\end{equation}
At the array, we have by definitions \eqref{eq:F1} and \eqref{eq:A3} that 
\begin{equation}
w(\vx) = \cN \phi(\vx), \qquad \vx \in \cA.
\label{eq:A9}
\end{equation}
Moreover, definition \eqref{eq:F11} and equation \eqref{eq:A8} give 
\begin{equation}
v(\vx) = \cS^{-1} w(\vx) = \cS^{-1} \cN \phi(\vx) = \cF \phi(\vx), \qquad 
\vx \in \cA.
\label{eq:A10}
\end{equation}
This proves that for $\phi$ given in \eqref{eq:A1}, we have 
\begin{equation}
\cF \phi(\vx) = \int_{\cA}dS_{\vy} \, \overline{\us(\vx,\vy)
  \phi(\vy)}, \qquad \vx \in \cA. \label{eq:A11}
\end{equation}

It remains to prove the result for functions 
\begin{equation*}
\phi(\vx) = \phi^{(1)}(\vx) + i \phi^{(2)}(\vx), \qquad \vx \in \cA,
\end{equation*}
 with $\phi^{(l)}$ defined in \eqref{eq:F16}. These have the same
 expression as \eqref{eq:A1}, so we write directly from \eqref{eq:A11}
 that
\begin{equation*}
\cF \phi^{(l)}(\vx) = \int_{\cA}dS_{\vy} \, \overline{\us(\vx,\vy)
  \phi^{(l)} (\vy)}, \qquad \vx \in \cA, ~~ l = 1,2. 
\end{equation*}
Because $\phi^{(l)}$ satisfy equation \eqref{eq:A2}, for $l = 1,2$,
  we have that
\begin{align*}
\int_{\cA} dS_{\vy} \,G(\vx,\vy) \big[\phi^{(1)}(\vy) + i
  \phi^{(2)}(\vy)\big] = \int_{\cA} dS_{\vy} \,\overline{G(\vx,\vy)
  \big[\phi^{(1)}(\vy) -i \phi^{(2)}(\vy)\big]}, \qquad \vx \in
\partial \Omega. 
\end{align*}
Then, the analogues of (\ref{eq:A3}--\ref{eq:A4}),
\begin{align*}
w(\vx) &= \int_{\cA} dS_{\vy} \, \us(\vx,\vy) \big[\phi^{(1)}(\vy) + i
  \phi^{(2)}(\vy)\big], \qquad \vx \in \cW \setminus
\overline{\Omega}, \\ v(\vx) &= \int_{\cA} dS_{\vy} \,
\overline{\us(\vx,\vy) \big[\phi^{(1)}(\vy) - i \phi^{(2)}(\vy)\big]},
\qquad \vx \in \cW \setminus \overline{\Omega},
\end{align*}
satisfy (\ref{eq:A5}--\ref{eq:A8}), and we conclude as above that 
\begin{equation*}
v(\vx) = \cS^{-1} w(\vx) = \cS^{-1} \cN \phi(\vx) = \cF \phi(\vx) =
\cF \phi^{(1)}(\vx) + i \cF\phi^{(2)}(\vx) , \qquad \vx \in \cA.
\end{equation*}
This proves Lemma \ref{lem.1}. $\Box$

\section{Proof of Theorem \ref{thm.1}}
\label{ap:B}
We begin in section \ref{ap:B.1} with the proof of \eqref{eq:F13}. 
The proofs of statements (i) and (ii) of the theorem are in sections
\ref{ap:B.2} and \ref{ap:B.3}. We use throughout the appendix the notation 
\[
\cW_{A}^+ = (x_\cA,0) \times \cX ~~ \mbox{and} ~~ 
\cW_{A}^- = (-\infty, x_\cA) \times \cX.
\]

\subsection{The factorization of $\cN$}
\label{ap:B.1}
Consider the operator $\cM:H^{\frac{1}{2}}(\partial \Omega)
\to L^2(\cA)$, defined by
\begin{equation}
\cM f(\vx) = w(\vx), \qquad \vx \in \cA, ~~ \forall \, f \in
H^{\frac{1}{2}}(\partial \Omega),
\label{eq:B2}
\end{equation}
where $w \in \mathscr{H}(W\setminus \overline{\Omega})$ satisfies 
the outgoing radiation condition at range $x < x_{\Omega}$ and the boundary 
condition 
\begin{equation}
  w(\vx) = - f(\vx), \qquad \vx \in \partial \Omega.
\label{eq:B3}
\end{equation}
By the definition of $\us(\vx,\vy)$ and using $f =G(\cdot, \vy)$ in
\eqref{eq:B3}, we have
\begin{equation*}
\us(\vx,\vy) = \big[\cM G(\cdot, \vy)|_{\partial \Omega}\big] (\vx),
\qquad \vx \in \cA.
\end{equation*}
Therefore, definitions (\ref{eq:F1}--\ref{eq:F2}) and the linearity of
$\cM$ give
\begin{align*}
\cN g(\vx) &= \int_{\cA} dS_{\vy} \big[\cM G(\cdot, \vy)|_{\partial
    \Omega}\big] (\vx) g(\vy) = \cM \Big[ \int_{\cA} dS_{\vy}
  G(\cdot, \vy)|_{\partial \Omega} g(\vy)\Big](\vx) \\ &= \cM \cT
g(\vx), \qquad \vx \in \cA, ~~ \forall \, g \in L^2(\cA).
\end{align*}
This proves the factorization 
\begin{equation}
\cN = \cM \cT.
\label{eq:B5}
\end{equation}

It remains to prove that 
\begin{equation}
{\cS^{-1} \cM = \cT^\star \Lambda.}
\label{eq:B6}
\end{equation}
Take any $f \in H^{\frac{1}{2}}(\partial \Omega)$ and use it to define
$h \in H^{-\frac{1}{2}}(\partial \Omega)$ by
\begin{equation*}
h(\vx) = \Lambda f(\vx), \qquad \vx \in \partial \Omega.
\end{equation*}
With this $h$, we obtain from definition \eqref{eq:F3} that
\begin{equation*}
\cT^\star h(\vx) = \cT^\star \Lambda f(\vx) = \int_{\partial \Omega} dS_{\vz} \, 
\overline{G(\vz,\vx)} h(\vz), \qquad \vx \in \cA.
\end{equation*}
If we let $v \in \mathscr{H}(W\setminus \overline{\Omega})$ be defined
by
\begin{equation}
v(\vx) = \int_{\partial \Omega} dS_{\vz} \, \overline{G(\vz,\vx)}
h(\vz), \qquad \vx \in \cW \setminus \overline{\Omega},
\label{eq:B9}
\end{equation}
then we have 
\begin{align}
v(\vx) &= \cT^\star \Lambda f(\vx), \qquad \vx \in \cA, \label{eq:B10}\\
v(\vx) &= -f(\vx), \qquad \vx \in \partial \Omega,
\label{eq:B11}
\end{align}
where \eqref{eq:B11} is obtained from definition \eqref{eq:F4}.  We
also have that $v|_{\cA} \in \mathscr{H}^{in}(\cA),$ so we can define
\begin{equation*}
{\cS^{-1}w(\vx) = v(\vx) = \cT^\star \Lambda f(\vx)}, \qquad \vx \in \cA.
\end{equation*} 
The factorization \eqref{eq:B6} follows from this equation, definition
 \eqref{eq:F11} of $\cS^{-1}$, equation \eqref{eq:B11} and definition
\eqref{eq:B2} of $\cM$, which give
\begin{equation*}
w(\vx) = \cM f(\vx), \qquad \vx \in \cA. 
\end{equation*}
$\Box$

\subsection{Proof of statement (i)}
\label{ap:B.2}
Consider the operator $\widetilde{\cT}:H^{-\frac{1}{2}}(\cA) \to
H^{\frac{1}{2}}(\partial \Omega)$, 
\begin{equation}
\widetilde{\cT} g(\vz) = \int_{\cA} d S_{\vy} \, G(\vz,\vy) g(\vy),
\qquad \forall \, \vz \in \partial \Omega, ~~ \forall \, g \in
L^2(\cA).
\label{eq:B19}
\end{equation}
whose restriction to the domain $L^2(\cA) \subset
H^{-\frac{1}{2}}(\cA)$ is the operator $\cT$ defined in \eqref{eq:F2}.

\subsubsection{Proof that $\widetilde\cT$ is bounded}
The Green's function $G(\vz,\vx)$ is smooth for $\vz \ne \vx$,
so  
\[
v(\vz) = \int_{\cA} dS_{\vx} \, G(\vz,\vx) g(\vx) \in
H^1\big(\cW_{\cA}^+\big), \qquad \forall \, g \in
H^{-\frac{1}{2}}(\cA).
\]
Moreover, for {$\vz \in \Omega$}, 
\[
\Delta_{\vz} v(\vz) = \int_{\cA} dS_{\vx} \, \Delta_{\vz} G(\vz,\vx)
g(\vx) = - k^2 \int_{\cA} dS_{\vx} \, G(\vz,\vx) g(\vx) = - k^2 v(\vz).
\]
By the mapping property of the single layer potential, $\|v\|_{H^1\big(\Omega\big)} \le C_1 \|g \|_{H^{-\frac{1}{2}}(\partial\Omega)}$ for some constant $C_ 1> 0$, which gives 
\[
\|\Delta_{\vz} v(\vz)\|_{L^2(\Omega)} \le k^2\| v(\vz)\|_{L^2(\Omega)} \le C_2  \|g \|_{H^{-\frac{1}{2}}(\partial\Omega)},
\]
for another constant $C_2 > 0$. Then, we obtain from 
\cite[Lemma 4.3]{mclean2000strongly}  that 
\begin{equation}
\widetilde{\cT} g = v|_{\partial \Omega} \in H^{\frac{1}{2}}(\partial
\Omega)~~ \mbox{and} ~ ~ \|\widetilde{\cT} g
\|_{H^{\frac{1}{2}}(\partial \Omega)} \le C
  \|g \|_{H^{-\frac{1}{2}}(\partial
\Omega)}, \label{eq:Bound}
\end{equation}
for yet another constant $C > 0$, so $\widetilde{\cT}$ is bounded.

\subsubsection{Proof that $\cT$ is compact}
Consider any bounded sequence in $L^2(\cA)$, which must have a weakly
convergent subsequence $\{g_n\}$ with weak limit $g \in L^2(\cA)$. 
Because $L^2(\cA)$ is compactly embedded in $H^{-\frac{1}{2}}(\cA)$, this 
sequence converges strongly in $H^{-\frac{1}{2}}(\cA)$, 
\begin{equation}
\lim_{n\to \infty} \|g_n -g\|_{H^{-\frac{1}{2}}(\cA)} = 0.
\label{eq:Blim}
\end{equation}
Recalling that $\cT$ is the restriction of $\widetilde{\cT}$ to the domain 
$L^2(\cA)$, we  have 
\[
\cT g_n = \widetilde{\cT} g_n ~~ \mbox{and}~~ \cT g = \widetilde{\cT} g, 
\]
and using \eqref{eq:Bound},
\begin{equation}
\|\cT g_n - \cT g \|_{H^{\frac{1}{2}}(\partial \Omega)} = \|\widetilde
  \cT g_n - \widetilde \cT g\|_{H^{\frac{1}{2}}(\partial \Omega)} \le
    C \|g_n -g\|_{H^{-\frac{1}{2}}(\cA)}.
\label{eq:Blim1}
\end{equation}
We conclude from (\ref{eq:Blim}--\ref{eq:Blim1}) that $\{\cT g_n\}$
converges to $\cT g$, strongly in $H^{\frac{1}{2}}(\partial
\Omega)$. This proves that $\cT$ is compact.

\subsubsection{Proof that $\cT$ is injective}
Let us define 
\begin{equation}
\label{eq:BI1}
v(\vz) = \int_{\cA} dS_{\vx}\, G(\vz,\vx) g(\vx), \qquad  \vz \in \cW,
\end{equation}
where $g$ satisfies 
\begin{equation}
v(\vz) = \cT g(\vz) = 0, \qquad \vz \in \partial \Omega.
\label{eq:BI1p}
\end{equation}
To prove injectivity, we must show that $g = 0$.

Equations (\ref{eq:BI1}--\ref{eq:BI1p}) give 
\begin{align*}
(\Delta_{\vz} + k^2) v(\vz) &= 0, \qquad \vz \in \Omega, \\
v(\vz) &= 0, \qquad \vz \in \partial \Omega,
\end{align*}
and by Assumption \ref{as.2}, 
\[
v(\vz) = 0, \qquad \forall \, \vz \in \overline{\Omega}.
\]
Since $v$ is analytic at $\vz \in (x_{\cA},0) \times \cX$, we obtain 
by unique continuation that 
\[
v(\vz) = 0,  \qquad \forall \, \vz \in \cW_{\cA}^+ \cup \cA.
\]
On the left of the array, we have 
\begin{align*}
(\Delta_{\vz} + k^2) v(\vz) &= 0, \qquad \vz \in \cW_{\cA}^-,
  \\ v(\vz) &= 0, \qquad \vz \in
  \cA,\\ \partial_{\vec{\boldsymbol{\nu}}_{\vz}} v(\vz) &=0, \qquad
  \vz \in \partial \cW,
\end{align*}
and $v$ satisfies the outgoing radiation condition at range $x <
x_{\cA}$. By the uniqueness of solution (see for example \cite[Lemma
  A.2]{borcea2018direct}) 
\[
v(\vz) = 0, \qquad \forall \, \vz \in \cW_{\cA}^-.
\]
But \eqref{eq:BI1} is a single layer potential, satisfying the jump
condition
\[
{-g|_{\cA}} = [\partial_{\vec{\boldsymbol{\nu}}} v ] \big|_{\cA} = 0,
\]
where $[\cdot ]\big|_{\cA}$ denotes the jump at $\cA$.
This proves that $\cT$ is injective. $\Box$

\subsection{Proof of statement (ii)}
\label{ap:B.3}
We show first that ${ -\mathscr{I}(\Lambda)}$ is positive semi-definite and then 
we prove the result on $\mathscr{R}(\Lambda)$.
\subsubsection{The operator { $\Im(\Lambda)$}}
Recall definition \eqref{eq:F4} and introduce  the functions
\begin{equation}
h(\vx) = \Lambda f(\vx) ~~\mbox{and} ~~ 
\widetilde h(\vx) = \Lambda \widetilde f(\vx), 
\qquad \vx \in \partial \Omega,
\label{eq:BII1}
\end{equation}
for arbitrary $f, \widetilde f \in H^{\frac{1}{2}}(\partial \Omega)$,
where $h \in H^{-\frac{1}{2}}(\partial \Omega)$ is the unique solution
of
\begin{equation}
\int_{\partial \Omega} dS_{\vy} \, \overline{G(\vx,\vy)}\, h(\vy) = -
f(\vx), \qquad \vx \in \partial \Omega,
\label{eq:BII2}
\end{equation}
and $\widetilde h$ satisfies a similar equation, with $\widetilde f$ in
the right hand side. 
Define 
\begin{align}
v(\vx) &= \int_{\partial \Omega} dS_{\vy} \, \overline{G(\vx,\vy)}\, 
h(\vy), \qquad \vx \in \cW \setminus \partial \Omega,
\label{eq:BII4}\\ \widetilde v(\vx) &= \int_{\partial \Omega} dS_{\vy} \,
\overline{G(\vx,\vy)} \, \widetilde h(\vy), \qquad \vx \in \cW
\setminus  \partial \Omega, \label{eq:BII5}
\end{align}
and note that \eqref{eq:BII2} implies
\begin{equation}
v(\vx) = - f(\vx) ~~\mbox{and} ~~ \widetilde v(\vx) = - \widetilde f(\vx), 
\qquad \vx \in \partial \Omega.
\label{eq:BR}
\end{equation}
Since (\ref{eq:BII4}--\ref{eq:BII5}) are single layer potentials, { we 
have from \cite[Theorem 6.11]{mclean2000strongly}}
\begin{equation}
{-h}\big|_{\partial \Omega} = [\partial_{\vec{\boldsymbol{\nu}}} v
]\big|_{\partial \Omega} ~~ \mbox{and}~~ { -\widetilde h}\big|_{\partial \Omega} =
[\partial_{\vec{\boldsymbol{\nu}}} \widetilde v ]\big|_{\partial \Omega},
\end{equation}
where $[\cdot ]$ denotes the jump at $\partial \Omega$.

These results imply that 
\begin{align} 
{-}\left< \Lambda f, \widetilde f \right>_{\partial \Omega} &= 
{ -}\left< h, \widetilde f \right>_{\partial \Omega} = { -}\int_{\partial \Omega} dS_{\vx} 
\, \overline{h(\vx)}\, \widetilde f(\vx) =  \int_{\partial \Omega} dS_{\vx} 
\, \overline{h(\vx)}\, \widetilde v(\vx), \nonumber \\
& =-
\int_{\partial \Omega} dS_{\vx} \, \big[ \overline{\partial_{\vnu}v^+(\vx)} -
  \overline{\partial_{\vnu}v^-(\vx)}\big] \widetilde v(\vx),
\label{eq:BII6}
\end{align}
with indexes $\pm$ denoting the function $v$ outside or inside
$\Omega$. Using the identity
\begin{align*}
\nabla_{\vx} \cdot \big[ \widetilde v(\vx) \nabla_{\vx}
  \overline{v^\pm(\vx)} \big] &= \widetilde v(\vx) \Delta_{\vx}
\overline{v^\pm(\vx)} + \nabla_{\vx} \widetilde v(\vx) \cdot
\nabla_{\vx} \overline{v^\pm(\vx)} \\&= - k^2 \widetilde
v(\vx) \overline{v^\pm(\vx)} + \nabla_{\vx} \widetilde v(\vx) \cdot
\nabla_{\vx} \overline{v^\pm(\vx)},
\end{align*}
and integration by parts, we obtain that 
\begin{align*}
\int_{\cW_{\cA}^+ \setminus \overline{\Omega}} d \vx \Big[ - k^2 \widetilde
  v(\vx) \overline{v^+(\vx)} + \nabla_{\vx} \widetilde v(\vx) \cdot
  \nabla_{\vx} \overline{v^+(\vx)}\Big] &= - \int_{\cA} dS_{\vx} \,
  \widetilde{v}(\vx) \overline{\partial_{\vnu} v^+(\vx)} \nonumber \\ &-
  \int_{\partial \Omega} dS_{\vx} \, \widetilde{v}(\vx)
  \overline{\partial_{\vnu} v^+(\vx)},
\end{align*}
and 
\begin{align*}
\int_{\Omega} d \vx \Big[ - k^2 \widetilde
  v(\vx) \overline{v^-(\vx)} + \nabla_{\vx} \widetilde v(\vx) \cdot
  \nabla_{\vx} \overline{v^-(\vx)} = &
  \int_{\partial \Omega} dS_{\vx} \, \widetilde{v}(\vx)
  \overline{\partial_{\vnu} v^-(\vx)}.
\end{align*}
Substituting these equations in \eqref{eq:BII6} we get 
\begin{align}
{-}\left< \Lambda f, \widetilde f \right>_{\partial \Omega} =& \int_{\cA}
dS_{\vx} \, \widetilde{v}(\vx) \overline{\partial_{\vnu} v(\vx)} +
\int_{\cW_{\cA}^+ \setminus \partial \Omega} d \vx \Big[ - k^2 \widetilde
  v(\vx) \overline{v(\vx)} + \nabla_{\vx} \widetilde v(\vx) \cdot
  \nabla_{\vx} \overline{v(\vx)}\Big], \label{eq:BIIp}
\end{align}
where we droped the $\pm$ indexes on $v$. The same calculation, 
with $v$ and $\widetilde v$ interchanged, gives 
\begin{align}
{ -}\left< \widetilde f, \Lambda^\star f \right>_{\partial \Omega} =&
{-}\left< \Lambda \widetilde f, f \right>_{\partial \Omega} = \int_{\cA}
dS_{\vx} \, {v}(\vx) \overline{\partial_{\vnu} \widetilde v(\vx)}
\nonumber \\ &+ \int_{\cW_{\cA}^+ \setminus \partial \Omega} \Big[ -
  k^2 v(\vx) \overline{\widetilde v^(\vx)} + \nabla_{\vx} v(\vx) \cdot
  \nabla_{\vx} \overline{\widetilde v(\vx)} = {-}\overline{ \left<
    \Lambda^\star f, \widetilde f \right>_{\partial
      \Omega}}.\label{eq:BIIpp}
\end{align}
Therefore, $\Im(\Lambda) =
\big(\Lambda-\Lambda^\star\big)/(2i)$ satisfies 
\begin{equation}
{-}\left< \Im(\Lambda) f, \widetilde f \right>_{\partial
  \Omega} = -\frac{1}{2 i} \int_{\cA} dS_{\vx} \Big[ \widetilde v(\vx)
  \overline{\partial_{\vnu} v(\vx)}- \overline{v(\vx)} \partial_{\vnu}
  \widetilde v(\vx)\Big].
\label{eq:BII9}
\end{equation}

We can write \eqref{eq:BII9} more explicitly  using the expression
\eqref{eq:uinc} of the Green's function in the definitions
(\ref{eq:BII4}--\ref{eq:BII5}) of $v$ and $\widetilde{v}$. We obtain 
\begin{align}
v(\vx) &= \sum_{j=0}^\infty \psi_j(\bxp) v_j \label{eq:BII10} \\
\partial_{\vnu} v(\vx) &= \sum_{j=0}^J i \beta_j \psi_j(\bxp) v_j + 
\sum_{j>J} |\beta_j| \psi_j(\bxp) v_j, \qquad \vx \in \cA, \label{eq:BII11}
\end{align}
where 
\begin{equation}
v_j = \left\{ \begin{array}{ll} 
-\frac{i}{\beta_j} \int_{\partial \Omega} dS_{\vy} \, h(\vy)
\psi_j(\byp) e^{i \beta_j x_{\cA}} \cos(\beta_j y), \qquad &j = 0, \ldots, J, \\\\
\frac{1}{|\beta_j|} \int_{\partial \Omega} dS_{\vy} \, h(\vy)
\psi_j(\byp) e^{|\beta_j| x_{\cA}} \cosh(|\beta_j| y), & j > J, 
\end{array} \right. 
\label{eq:BII12}
\end{equation}
and similar for $\widetilde v$. Substituting in \eqref{eq:BII11} and
using the orthogonality relation \eqref{eq:EigFOrt}, 
\begin{equation}
{-}\left< \Im(\Lambda) f, \widetilde f \right>_{\partial
  \Omega} = \sum_{j=0}^J \beta_j \widetilde v_j \overline{v_j}.
\label{eq:BII13}
\end{equation}
In particular, for $\widetilde f = f$, 
\begin{equation}
{-}\left< \Im(\Lambda) f, f \right>_{\partial \Omega} =
\sum_{j=0}^J \beta_j |v_j|^2 \ge 0, \qquad \forall \, f \in
H^{\frac{1}{2}}(\partial \Omega).
\label{eq:BII14}
\end{equation}
This proves that $-\Im(\Lambda)$ is  positive semi-definite. $\Box$

\subsubsection{The operator $\mathscr{R}(\Lambda)$} 
We introduce the operator $\Lambda_i$ by
\begin{equation}
h(\vx) = \Lambda_i f(\vx), 
\qquad \vx \in \partial \Omega,
\label{eq:BRR1}
\end{equation}
for arbitrary $f\in H^{\frac{1}{2}}(\partial \Omega)$, where $h \in H^{-\frac{1}{2}}(\partial \Omega)$ is the unique solution of
\begin{equation}
\int_{\partial \Omega} dS_{\vy} \, \overline{G_i(\vx,\vy)}\, h(\vy) = -
f(\vx), \qquad \vx \in \partial \Omega,
\label{eq:BRR2}
\end{equation}
with $G_i(\vx,\vy)$ the Green function when $k=i$. We let  $v_i$ satisfy (\ref{eq:BII4}) with $G(\vx,\vy)$ replaced by $G_i(\vx,\vy)$. By Assumption \ref{as.2}, both $\Lambda$ and $\Lambda_i$ have bounded inverses, and from \eqref{eq:F4} -- \eqref{eq:F5} we see that for any $h \in H^{-\frac{1}{2}}(\partial \Omega)$,
\begin{equation}
\Lambda^{-1} h (\vx)= - \int_{\partial \Omega} dS_{\vy} \, \overline{G(\vx,\vy)}\, h(\vy) , \qquad \vx \in \partial \Omega,
\label{eq:BRR3}
\end{equation}
The analogue of this equation holds for $\Lambda_i^{-1}$ with $G(\vx,\vy)$ replaced by $G_i(\vx,\vy)$.

Note that $G(\vx,\vy) - G_i(\vx,\vy)$ satisfies the Helmoltz equation and it is smooth. In particular, this is so for $\vx,\vy \in \partial \Omega$. Because $\overline{G(\vx,\vy)} - \overline{G_i(\vx,\vy)}$ is the kernel of operator $\Lambda^{-1} - \Lambda_i^{-1}$, it follows that  $\Lambda^{-1} - \Lambda_i^{-1}$ is compact. Furthermore 
\[\Lambda - \Lambda_i = \Lambda_i (\Lambda_i^{-1} - \Lambda^{-1}) \Lambda
\]
 is compact and $\mathscr{R}(\Lambda- \Lambda_i)$ is compact. The representation (\ref{eq:BIIp}) of $\Lambda_i$ (where we replace $k$ by $i$) gives $\Lambda_i^*=\Lambda_i$ and 
\begin{align*}
-\left< \Lambda_i f, f \right>_{\partial \Omega} &\ge \|v_i\|_{H^1(\cW_{\cA}^+ \setminus \partial
  \Omega)}^2 \ge C \|f\|_{H^{\frac{1}{2}}(\partial \Omega)}, \qquad \forall \, f \in
H^{\frac{1}{2}}(\partial \Omega).
\end{align*}
This yields that $-\mathscr{R}(\Lambda) =-\Lambda_i - \mathscr{R}(\Lambda- \Lambda_i)$ is the sum of 
a positive definite, self-adjoint operator and a compact operator.
$\Box$

\section{Proof of Lemma \ref{lem.2}}
\label{ap:C}

Let us start with the case $\vz \in \Omega$. Since
$G(\cdot,\vz)\big|_{\partial \Omega}$ is in $H^{\frac{1}{2}}(\partial
\Omega)$, we can define $h \in H^{-\frac{1}{2}}(\partial \Omega)$ by
\begin{equation}
h \big|_{\partial \Omega} = - \Lambda \, \overline{G(\cdot,\vz)}\big|_{\partial \Omega},
\label{eq:C1}
\end{equation}
where we recall from definition \eqref{eq:F4} that $h$ is the unique
solution of
\begin{equation}
\int_{\partial \Omega} dS_{\vy} \, G(\vx,\vy) \overline{h(\vy)} =
G(\vx,\vz), \qquad \vx \in \partial \Omega.
\label{eq:C2}
\end{equation}
With this $h$, let 
\begin{equation}
w(\vx) = \int_{\partial \Omega} dS_{\vy} \, G(\vx,\vy)
\overline{h(\vy)}, \qquad \vx \in \cW \setminus \partial \Omega.
\label{eq:C3}
\end{equation}
Then, $v(\vx) = w(\vx) - {G(\vx,\vz)}$ is in $\mathscr{H}(\cW \setminus
\overline{\Omega})$ and it satsifies the outgoing radiation condition
at range $x < x_{\Omega}$ and the boundary condition $v\big|_{\partial
  \Omega} = 0$. By Assumption \ref{as.1}, we conclude that $
v(\vx) = 0 $ in $ \cW \setminus \overline{\Omega}.$ This implies in
particular that
\begin{equation}
w(\vx) = G(\vx,\vz), \qquad \vx \in \cA.
 \label{eq:C5}
\end{equation}
Furthermore, by definition \eqref{eq:F3}, we get for all $\vx \in \cA$, 
\begin{align}
\cT^\star h(\vx) = \int_{\partial \Omega} dS_{\vy} \, 
\overline{G(\vy,\vx)} h(\vy) = \int_{\partial \Omega} dS_{\vy} \, 
\overline{G(\vx,\vy)} h(\vy)  = \overline{w(\vx)} = \overline{G(\vx,\vz)}, 
\end{align}
where we used \eqref{eq:C3}, \eqref{eq:C5} and
the reciprocity of the Green's function.  This shows that
$\overline{G(\cdot, \vz)}\big|_{\cA} \in \mbox{range}(\cT^\star)$.

To prove the converse, suppose that $\vz \notin \Omega$ and assume for
a contradiction argument that $ \overline{G(\cdot, \vz)}\big|_{\cA} \in
\mbox{range}(\cT^\star).$ Then, there exists $h \in
H^{-\frac{1}{2}}(\partial \Omega)$ such that
\begin{equation}
\cT^\star h(\vx) = \int_{\partial \Omega} dS_{\vy} \,
\overline{G(\vy,\vx)} h(\vy) = \overline{G(\vx,\vz)}, \qquad \vx \in
\cA.
\label{eq:C6}
\end{equation}
This $h$ defines a function $w$ as in \eqref{eq:C3}, satisfying $w \in
H_{\rm loc}^1(\cW \setminus \overline \Omega)$, with trace
$w\big|_{\partial \Omega} \in H^{\frac{1}{2}}(\partial \Omega)$. If we
define further
\begin{equation}
v(\vx) = w(\vx) - G(\vx,\vz),
\label{eq:C7}
\end{equation}
then we obtain that it satisfies the boundary value problem
\begin{align*}
(\Delta_{\vx} + k^2) v(\vx) &= 0, \qquad \vx \in \cW_{A}^-, \\
\partial_{\vnu} v(\vx) &=0, \qquad \vx \in \partial \cW, \\
v(\vx)&=0, \qquad \vx \in \cA, 
\end{align*}
and the outgoing radiation condition at $x< x_{\cA}.$ This problem has
the unique solution (see for example \cite[Lemma
  A.2]{borcea2018direct})
\[
v(\vx) = 0, \qquad \vx \in \cW_{\cA}^- \cup \cA,
\]
and since $v$ is analytic at $\vx \notin \overline{\Omega} \cup
\{\vz\}$, we have by unique continuation 
\[
v(\vx) = 0 ~~ \mbox{i.e.}~~ w(\vx) = G(\vx,\vz), \qquad \vx \in \cW
\setminus \{ \overline{\Omega} \cup \{\vz\}\}.
\]
Then, $w(\vx)$ blows up, like $G(\vx,\vz)$, as $\vx \to
\vz$. This contradicts that
$w \in H_{\rm loc}^1(\cW \setminus \overline{\Omega})$ and
$w\big|_{\partial \Omega} \in H^{\frac{1}{2}}(\partial
\Omega)$. Therefore, $\overline{G(\cdot, \vz)}\big|_{\cA} \notin
\mbox{range}(\cT^\star)$ when $\vz \notin \Omega$. $\Box$

\section{Proof of Theorem \ref{thm.NF}}
\label{ap:D}
To prove the theorem, we begin with two lemmas, proved in sections 
\ref{ap:D1} and \ref{ap:D2}.

\vspace{0.05in}
\begin{lemma}
\label{lem.D1} 
Denote by $\overline{\cT(\csP)}$ the closure of the image of the set
$\csP$ defined in \eqref{eq:F26} under the operator $\cT$ defined in
\eqref{eq:F2}. Recall also the sets $\csP_0$ and $\csP_0^\perp$
defined in \eqref{eq:R4} and denote by $\cT(\csP_0)$ and
$\cT(\csP_0^\perp)$ their image under $\cT$. We have
\begin{equation}
\left< -\Im(\Lambda) f,f \right>_{\partial \Omega} \ne 0, \qquad 
\forall \, f \in \overline{\cT(\csP)}, ~~ f \notin \cT(\csP_0), ~~ f \ne 0.
\label{eq:D1}
\end{equation}
Moreover, there exists a positive constant $C$ such that
\begin{equation}
\left< - \Im(\Lambda) f,f \right>_{\partial \Omega} \ge C
\|f\|^2_{H^{\frac{1}{2}}(\partial \Omega)}, \qquad \forall \, f \in
\cT(\csP_0^\perp), ~~ f \ne 0.
\label{eq:D2}
\end{equation}
\end{lemma}

\vspace{0.05in}
\begin{lemma}
\label{lem.D2} 
A search point $\vz$ lies in $\Omega$ and therefore, by Lemma \ref{lem.2},
we have $\overline{G(\cdot, \vz)}\big|_{\cA} \in \mbox{range}(\cT^\star)$, if
and only if $\overline{G_{\csP}(\cdot, \vz)}\big|_{\cA} \in
\mbox{range}(\cT^\star)\big|_{\csP}.$
\end{lemma}

\vspace{0.07in} 
Proof of  Theorem \ref{thm.NF}: 
By assumption, $\overline{G_{\csP}(\cdot, \vz)}\big|_{\cA} \notin
\csP_0$, so there exists $\varphi \in \csP_{0}^\perp$ so that 
\begin{equation*}
\big( \overline{G_{\csP}(\cdot, \vz)}, \varphi \big)_{\cA} \ne 0.
\end{equation*}
Therefore,
\begin{equation*}
\Phi = \{\varphi \in \csP_{0}^\perp, ~~ \big( \overline{G_{\csP}(\cdot,
    \vz)}, \varphi \big)_{\cA} = 1 \} \ne \emptyset.
\end{equation*}
Because $\vz \in \Omega$, we conclude from Lemma \ref{lem.D2} that
$\overline{G_{\csP}(\cdot, \vz)}\big|_{\cA} \in
\mbox{range}(\cT^\star)\big|_{\csP}.$ That is to say,
\begin{equation*}
\exists \, \theta \in H^{-\frac{1}{2}}(\partial \Omega) ~~ \mbox{such that}~~ 
\overline{G_{\csP}(\vx,
  \vz)} = \PP \cT^\star \theta (\vx), \qquad \vx \in \cA.
\end{equation*}
We also have from the factorization of $\cF$ in Theorem \ref{thm.1}
and definition \eqref{eq:F13p} that 
\begin{equation*}
 { -\big(\Im(\cF) \varphi,\varphi\big)_{\cA} = -\big(\cT^\star
    \Im(\Lambda) \cT \varphi,\varphi\big)_{\cA} = -\left<
    \Im(\Lambda) \cT \varphi,\cT \varphi \right>_{\partial
      \Omega} }\ge C \|\cT \varphi \|^2_{H^{\frac{1}{2}}(\partial
      \Omega)}, 
\end{equation*}
for all $\varphi \in \Phi$, where we used the bound \eqref{eq:D2} in
Lemma \ref{lem.1}. With these results we get 
\begin{align*}
C &= C \Big| \big(\overline{G_{\csP}(\cdot, \vz)},\varphi\big)_{\cA}
\Big|^2 = C \Big| \big(\PP \cT^\star \theta,\varphi \big)_{\cA}\Big|^2
= C \Big| \big(\cT^\star \theta,\varphi \big)_{\cA}\Big|^2 \nonumber
\\ &= C \Big| \left< \theta,\cT \varphi \right>_{\partial
  \Omega}\Big|^2 \le C
\|\theta\|^2_{H^{-\frac{1}{2}}(\partial \Omega)} \|\cT \varphi
\|^2_{H^{\frac{1}{2}}(\partial \Omega)} \nonumber \\ &\le
\|\theta\|^2_{H^{-\frac{1}{2}}(\partial \Omega)} \big( { -\Im(\cF)}
\varphi,\varphi\big)_{\cA}, \qquad \forall \, \varphi \in \Phi.
\end{align*}
and \eqref{eq:R6} follows. $\Box$

\subsection{Proof of Lemma \ref{lem.D1}}
\label{ap:D1}
To prove statement \eqref{eq:D1}, we use a contradiction argument.
Suppose that
\begin{equation}
\exists \, f \in \overline{\cT(\csP)} \setminus \cT(\csP_0), ~~ f \ne 0,
~~ \mbox{such that} ~~ \left< \Im(\Lambda) f,f \right>_{\partial
  \Omega} = 0.
\label{eq:D1.1}
\end{equation}
With this $f$, we define
\begin{equation}
h(\vx) = \Lambda f(\vx), \qquad \vx \in \partial \Omega,
\label{eq:D1.2}
\end{equation}
where we recall from definition \eqref{eq:F4} that $h \in
H^{-\frac{1}{2}}(\partial \Omega)$ is the unique solution of 
\begin{equation}
\int_{\partial \Omega} dS_{\vy} \, \overline{G(\vx,\vy)} h(\vy) = -
f(\vx), \qquad \vx \in \partial \Omega.
\label{eq:D1.3}
\end{equation}
Define also 
\begin{equation}
w(\vx) = \int_{\partial \Omega} dS_{\vy} \, G(\vx,\vy)
\overline{h(\vy)}, \qquad \vx \in \cW \setminus \partial \Omega,
\label{eq:D1.4}
\end{equation}
and note that it is like the complex conjugate of \eqref{eq:BII4}. Then, 
\eqref{eq:BII14} gives 
\[
{-\left< {\mathscr{I}}(\Lambda) f, f \right>_{\partial \Omega} }=
\sum_{j=0}^J \beta_j |w_j|^2, ~~ ~\mbox{with}~~ w_j = \int_{\cA}
dS_{\vx} \, w(\vx) \psi_j(\bxp), \qquad j = 0, \ldots, J,
\]
and assumption \eqref{eq:D1.1} implies that $w\big|_{\cA}$ is purely
evanescent. Therefore, using definitions (\ref{eq:F6}--\ref{eq:F8}),
we have
\begin{equation}
w, \overline{w} \in \mathscr{H}(W\setminus \overline{\Omega}) ~~
\mbox{and} ~~ w\big|_{\cA}, \overline{w}\big|_{\cA} \in
\mathscr{H}^{\rm out}(A).
\label{eq:D1.5}
\end{equation}
Moreover, equation \eqref{eq:D1.3} gives
\begin{equation}
\overline{w(\vx)} = - f(\vx), \qquad \vx \in \partial \Omega.
\label{eq:D1.6}
\end{equation}

Since $f \in \overline{\cT(\csP)}$, there is a sequence $\{g_n\}$ in $\csP$ such that the sequence $\{f_n\}$ defined by
\[
f_n(\vx) = \cT g_n(\vx), \qquad \vx \in \partial \Omega,
\]
converges to $f$. The convergent sequence $\{f_n\}$ must be bounded.  Because
$\cT$ is linear and injective, $\cT: \csP \to \cT(\csP)$ is invertible
and the inverse $\cT^{-1}:\cT(\csP) \to \csP$ is also a linear
operator. Moreover, since $\csP$ is finite dimensional, so is
$\cT(\csP)$. Thus, $\cT^{-1}$ is a map between finite dimensional
spaces, which means that it can be represented by a matrix and it is
bounded. We conclude that the sequence $\{g_n\}$, with $g_n
= \cT^{-1} f_n$ is bounded. Then, by the Bolzano-Weierstrass theorem,
there is a subsequence, still denoted by $\{g_n\}$ that converges to
$g \in \csP$, and we must have 
\begin{align}
f(\vx) = \cT g(\vx) = \int_{\cA} dS_{\vy} \, G(\vx,\vy) g(\vy) = -
\int_{\cA} dS_{\vy} \, \us(\vx,\vy) g(\vy), \qquad \vx \in \partial
\Omega.
\label{eq:D1.7p}
\end{align}
Here we used definition \eqref{eq:F2} and equation \eqref{eq:usc3}.

Note that for $\vx \in \cW \setminus \overline{\Omega}$, 
\begin{equation}
\int_{\cA} dS_{\vy} \, \us(\vx,\vy) g(\vy) \in \mathscr{H}(W\setminus
\overline{\Omega}), 
\label{eq:D1.7}
\end{equation}
and for $\vx \in \cA$,
\begin{equation}
\int_{\cA} dS_{\vy} \, \us(\vx,\vy) g(\vy) \in \mathscr{H}^{\rm
  out}(A).
\label{eq:D1.8}
\end{equation}
Equations (\ref{eq:D1.5}--\ref{eq:D1.6}) and
(\ref{eq:D1.7p}--\ref{eq:D1.8}) and the uniqueness of solutions imply 
\[
\overline{w(\vx)} = \int_{\cA} dS_{\vy} \, \us(\vx,\vy) g(\vy) =
\int_{\partial \Omega} dS_{\vy} \, \overline{G(\vx,\vy)} h(\vy),
\qquad \vx \in \cW \setminus \overline{\Omega}.
\]
However, we concluded above that $w\big|_{\cA}$ is purely evanescent,
which means that
\[
\PP \cN g(\vx) = \PP \int_{\cA} dS_{\vy} \, \us(\vx,\vy) g(\vy) = 0,
\qquad \vx \in \cA.
\]
This contradicts that $f = \cT g \notin \cT(\csP_0)$, and completes
the proof of \eqref{eq:D1}.

To prove statement \eqref{eq:D2}, we also argue by contradiction. Let
us work with the normalized functions \[ \varphi =
{f}/{\|f\|_{H^{\frac{1}{2}}(\partial \Omega)}}.  \] If \eqref{eq:D2} is
not true, then for any $n \in \mathbb{N}$, we can find $\varphi_n \in
\cT(\csP_0^\perp)$ with norm $\|\varphi_n\|_{H^{\frac{1}{2}}(\partial
  \Omega)} = 1$ such that
\begin{equation}
0 \le \left< { -\Im(\Lambda)} \varphi_n, \varphi_n \right>_{\partial
  \Omega} < \frac{1}{n}.
\label{eq:D1.12}
\end{equation}
Because $\varphi_n \in \cT(\csP_0^\perp)$, we can define a new
sequence $\{g_n\}$ in $\csP_0^\perp$,
\[
g_n = T^{-1} \varphi_n, \qquad \forall \, n \in \mathbb{N},
\]
which is bounded because $\cT^{-1}:\cT(\csP) \to \csP$ is bounded.
Then, by the Bolzano-Weierstrass theorem there is a subsequence, still
denoted by $\{g_n\}$, which converges to $g \in \csP_0^\perp$.  This
$g$ cannot be zero because $ \varphi = \cT g$ is the limit of the
sequence $\{\varphi_n\}$ of norm one. Taking the $n \to \infty $ limit
in \eqref{eq:D1.12} we get 
\[
\left< \Im(\Lambda) \varphi, \varphi \right>_{\partial \Omega} = 0,
\]
which contradicts statement \eqref{eq:D1}. Thus, statement
\eqref{eq:D2} must be true. $\Box$

\subsection{Proof of Lemma \ref{lem.D2}}
\label{ap:D2}
If $\overline{G(\cdot,\vz)}\big|_{\cA} \in \mbox{range}(\cT^\star)$ it
is obvious, from definitions, that
$\overline{G_{\csP}(\cdot,\vz)}\big|_{\cA} \in
\mbox{range}(\cT^\star)\big|_{\csP}.$ Thus, let us prove the converse. 

For a proof by contradiction, suppose that
$\overline{G_{\csP}(\cdot,\vz)}\big|_{\cA} \in
\mbox{range}(\cT^\star)\big|_{\csP}$ and yet,
\begin{equation}
\overline{G(\cdot,\vz)}\big|_{\cA} \notin \mbox{range}(\cT^\star).
\label{eq:D2.1}
\end{equation}
This means, by Lemma \ref{lem.2} that $\vz \notin \Omega$.
Then, there is $h \in H^{-\frac{1}{2}}(\partial \Omega)$ satisfying 
\begin{equation}
\PP \cT^\star h(\vx) = \overline{G_{\csP}(\vx,\vz)}, \qquad \vx \in \cA.
\label{eq:D2.2}
\end{equation}
With this $h$, we define 
\begin{align}
w(\vx)&=\int_{\partial \Omega} dS_{\vy} \, G(\vx,\vy)
\overline{h(\vy)}, \qquad \vx \in \cW \setminus \overline{\Omega},
\\ w_{\csP}(\vx)&=\int_{\partial \Omega} dS_{\vy}\, G_{\csP}(\vx,\vy)
\overline{h(\vy)}, \qquad \vx \in \cW \setminus \overline{\Omega},
\label{eq:D2.3p}
\end{align}
and obtain from \eqref{eq:D2.2} and definition \eqref{eq:F3} that
\begin{equation}
w_{\csP}(\vx) = \PP w(\vx) = \overline{\PP \cT^\star h(\vx)} =
G_{\csP}(\vx,\vz), \qquad \vx \in \cA.
\end{equation}
Note that $G_{\csP}(\vx,\vz)$ and $w_{\csP}$ solve the same problem in
$\cW_{A}^-$, with the same outgoing radiation condition. By the uniqueness
of solutions, we must have 
\begin{equation}
w_{\csP}(\vx) = G_{\csP}(\vx,\vz), \qquad \vx \in \cW_{\cA}^{-}.
\label{eq:D2.3}
\end{equation}
On the right of the array, at $\vx \notin \overline{\Omega}
\cup\{\vz\}$, $G_{\csP}(\vx,\vz)$ and $w_{\csP}$ again solve the same
problem, so by unique continuation of \eqref{eq:D2.3} we have
\begin{equation}
w_{\csP}(\vx) = G_{\csP}(\vx,\vz), \qquad \vx \in \cW \setminus
\{\overline{\Omega} \cup\{\vz\}\}.
\label{eq:D2.4}
\end{equation}
However, definition \eqref{eq:D2.3p} implies that $w_{\csP}$ and
$\partial_{x}^2 w_{\csP}$ are smooth in $\cW \setminus
\overline{\Omega}$, whereas
\[
\partial_x^2 G_{\csP}(\vx,\vz) = - \sum_{j=0}^J \psi_j(\bxp)
\psi_j(\bzp) \Big[ \delta(x-z) + i \beta_j \Big(e^{i \beta_j|x-z|} +
  e^{i \beta_j|x+z|}\Big) \Big]
\]
has a Dirac delta singularity at $\vz \in \cW \setminus
\overline{\Omega}$ with range $z = x$. We reached a contradiction, so
\eqref{eq:D2.1} cannot be true. $\Box$.

\bibliographystyle{plain} \bibliography{PAPER}

\end{document}